\documentclass[a4paper,12pt,reqno]{amsart}

\usepackage{amsmath,amssymb,amsthm}
\usepackage{latexsym}
\usepackage{color}
\usepackage{graphicx}
\usepackage{mathrsfs}
\usepackage{enumerate}
\usepackage[abbrev]{amsrefs}
\usepackage[T1]{fontenc}

\usepackage{mathtools}
\mathtoolsset{showonlyrefs=true}

\allowdisplaybreaks[4]

\theoremstyle{plain}
\newtheorem{thm}{Theorem}[section]
\newtheorem{lemm}[thm]{Lemma}
\newtheorem{prop}[thm]{Proposition}
\newtheorem{cor}[thm]{Corollary}

\theoremstyle{definition}
\newtheorem{df}[thm]{Definition}
\newtheorem{rem}[thm]{Remark}

\makeatletter

\@addtoreset{equation}{section}
\makeatother

\renewcommand{\div}{\operatorname{div}}

\renewcommand{\leq}{\leqslant}
\renewcommand{\geq}{\geqslant}

\newcommand{\n}[1]{{\left\|#1\right\|}}
\newcommand{\lp}[1]{\left[#1\right]}
\newcommand{\Mp}[1]{\left\{#1\right\}}

\newcommand{\R}{\mathbb{R}}

\newcommand{\Z}{\mathbb{Z}}
\newcommand{\ep}{\mathbb{\varepsilon}}
\newcommand{\wt}{\widetilde}
\newcommand{\Dh}{\dot{\Delta}}

\newcommand{\Jx}{\langle x_3 \rangle}
\newcommand{\Jxx}{\langle x \rangle}
\newcommand{\Jy}{\langle y_3\rangle}
\newcommand{\Jxy}{\langle x_3-y_3\rangle}
\newcommand{\JX}{\langle x \rangle}

\begin{document}
\title[3D stationary Navier--Stokes equations]
{Decay rates of three dimensional stationary Navier--Stokes flows at the spatial infinity}
\author[M.~Fujii]{Mikihiro Fujii}
\address[M.~Fujii]{Graduate School of Science, Nagoya City University, Nagoya, 467-8501, Japan}
\email[M.~Fujii]{fujii.mikihiro@nsc.nagoya-cu.ac.jp}
\author[H.~Tsurumi]{Hiroyuki Tsurumi}
\address[H.~Tsurumi]{Graduate School of Technology, Industrial and Social Sciences, Tokushima University, Tokushima 770--8506, Japan}
\email[H.~Tsurumi]{tsurumi.hiroyuki@tokushima-u.ac.jp}
\author[X.~Zhang]{Xin Zhang}
\address[X.~Zhang]{School of Mathematical Sciences, 
Key Laboratory of Intelligent Computing and Applications (Ministry of Education), Tongji University, No.1239, Siping Road, Shanghai (200092), China}
\email[X.~Zhang]{xinzhang2020@tongji.edu.cn}
\keywords{{ the three dimensional stationary Navier--Stokes equations; 
well-posedness; the critical hybrid Besov spaces; pointwisely decay estimates.}
}
\subjclass[2020]{35Q35, 76D03, 76D05}
\begin{abstract}
In this paper, we establish the well-posedness results of the three dimensional stationary Navier--Stokes equations (SNS) in some \emph{critical} hybrid type Besov spaces with respect to the scaling invariant structure of (SNS).
Although such critical functional spaces contain the functions with singularities, we give some sufficient conditions such that the $L^{\infty}$-norm of the solutions of (SNS) decay at the infinity within some polynomial type rate.
\end{abstract}
\maketitle


\section{Introduction}
\subsection{Model, previous works and motivation}
In this paper, we consider the three dimensional stationary incompressible Navier--Stokes equations on the whole space $\mathbb{R}^3$:
\begin{align}
\label{eq:NS}
    \begin{cases}
        -\Delta u + (u \cdot \nabla)u + \nabla p = \div f, \qquad & x \in \mathbb{R}^3,\\
        \div u = 0, \qquad & x \in \mathbb{R}^3.
    \end{cases}
\end{align}
Here the unknowns $u=(u_1(x),u_2(x),u_3(x)) $ and $p=p(x)$ denote the velocity vector and the pressure of the fluid at the point $x=(x_1,x_2,x_3)\in \mathbb{R}^3$ respectively. Moreover, $\div f$ denotes the external force of the divergence form, where $f=(f_{j,k}(x))_{1\leq j,k \leq 3}$ is a given $\mathbb{R}^{3 \times 3}$-valued function.

For a long time, many research results have been presented regarding the existence, uniqueness, and regularity of solutions to the system \eqref{eq:NS}.
Classically, Leray \cites{Ler-33}, Ladyzhenskaya \cites{Lad-59} and Fujita \cites{Fuj-61} proved the existence of solutions, and Finn \cites{Fin-65} investigated solutions in an exterior domain.
On the other hand, Heywood \cites{Hey-70} proved that the non-stationary solution $u(t)$ of Navier--Stokes system converges to the solution $u$ of \eqref{eq:NS} as the time $t\to\infty$.

As the theories on functional spaces and harmonic analysis has developed, the solutions of \eqref{eq:NS} has been considered in the framework of the so-called \emph{critical} spaces.
To be precise, let us denote 
\begin{equation}\label{scaling}
u\leadsto u_\lambda:= \lambda u( \lambda x),\quad
p\leadsto p_\lambda:= \lambda^2 p ( \lambda x),\quad 
f\leadsto f_\lambda:=\lambda^2 f (\lambda x)
\end{equation} 
for any constant $\lambda>0.$
Clearly, if the triplet $(u,p,f)$ satisfies \eqref{eq:NS}, then so does $(u_\lambda, p_\lambda, f_\lambda)$.
Thus it is very natural to seek for the solutions of \eqref{eq:NS} in some functional spaces which are invariant under the scaling transformations in \eqref{scaling}. In such situation, we say such functional spaces are critical for the problem \eqref{eq:NS} and the resulting solutions are critical as well. 
To the present, various researchers have established the well-posedness theory of \eqref{eq:NS} in different critical spaces,
such as Secchi \cites{Sec-88} and Chen \cites{Che-93} in the Lebesgue space $L^3(\R^3),$ 
Kozono--Yamazaki \cites{Koz-Yam-95-PJA, Koz-Yam-95-IUMJ} and Bjorland--Brandolese--Iftimie--Schonbek \cites{Bjo-Bra-Ift-Sch-11} in the Morrey space $\dot{M}_{3,r}$, and Kaneko--Kozono--Shimizu \cites{Kan-Koz-Shi-19} and T. \cites{Tsu-19-DIE} in the standard homogenous Besov space $\dot{B}^{3/p-1}_{p,q}(\R^3)$ and Triebel--Lizorkin space $\dot{F}^{3/p-1}_{p,q}(\R^3)$, respectively.
Furthermore, the well-posedness issue here means the uniquely existence of small solution with respect to any given small force. 
More precisely, we say the system \eqref{eq:NS} is well-posed from the space $D$ of the datum $f$ to the space $S$ of solutions $u$ if the following statement holds: 
There exist constants $\varepsilon>0$ and $\delta=\delta(\varepsilon)>0$ such that the solution map
\footnote{For any Banach space $(X,\|\cdot\|_X)$ and the constant $r>0$, $B_X(r):=\{x\in X; \|x\|_X<r\}$ denotes the open ball with radius $r$ in $X$.} of \eqref{eq:NS}
$$f\in B_D(\varepsilon) \mapsto u\in B_S(\delta)$$  
is well-defined and continuous with regard to those norms.
For the ill-posedness results in critical Besov spaces with too weak norms, we refer to T. \cites{Tsu-19-JMAA, Tsu-19-ARMA}, Li--Yu--Zhu \cites{Li-Yu-Zhu} and Tan--T.--Z. \cites{Tan-Tsu-Zha} for $\R^3$. 
On the other hand, it was also investigated by F. \cites{Fujii-APDE} that in the whole plane $\R^2$, \eqref{eq:NS} is ill-posedness in the critical Besov space $\dot{B}^{2/p-1}_{p,1}(\R^2)$ with $1\leq p \leq 2$.

However, despite of many studies on the well-posedness, there seems no previous results on the decay rate for the solutions of \eqref{eq:NS} at the spatial infinity. For the decay property of the solution, although the $L^3$-solutions in \cite{Che-93} may decay at the infinity in some sense, there is no detailed information on the decay rate of $|u(x)|$ as $|x|\to\infty,$ to say nothing of the solutions in \cites{Kan-Koz-Shi-19} with less regularities.

In order to construct the solutions of \eqref{eq:NS} with the (pointwise) decay property, we consider the well-posedness issue of the system \eqref{eq:NS} in the framework of the critical hybrid Besov spaces (see Definition \ref{def:Besov} for more details). 
Such hybrid Besov spaces can be regarded as the Chemin--Lerner type spaces (see \cite{Che-Ler-95} for instance) used in the non-stationary problem, where the horizontal direction $x_h:=(x_1,x_2)$ is estimated in the homogeneous Besov spaces and the vertical variable $x_3$ is bounded in the Lebesgue norms.
Thus, compared to the standard Besov spaces, it makes us possible to identify the spatial behavior of functions in the vertical direction by regarding $x_3$ as a ``time'' variable inspired by the work \cite{Fuj-pre}.

\subsection{Main results}
As in the standard analysis of the stationary Navier--Stokes problem, we apply the so-called \emph{Helmholtz} (or \emph{Leray})  \emph{projection}
$$\mathbb{P}=\left\{ \delta_{jk}+\partial_{x_j}\partial_{x_k}(-\Delta)^{-1}\right\}_{1\leq j,k\leq 3}$$
to \eqref{eq:NS}, and then we see that
\begin{align}\label{eq:NS2}
        -\Delta u = \mathbb{P}\div (f - u \otimes u), \qquad & x \in \mathbb{R}^3.
\end{align}
For any functional space $X$ of $\R^3$-valued mappings, we denote 
\begin{equation*}
    \mathbb{P}X:=\{ v\in X:  \div v =0\}.
\end{equation*}

The first result on the well-posedness issue of \eqref{eq:NS2} is as follows
\footnote{See Definition \ref{def:Besov} for the functional spaces used in our results.}:
\begin{thm}\label{thm:wp_1}
Let $(p,r)\in [1,\infty)^2$ and $q\in [1,\infty]$ satisfying the conditions 
\begin{equation*}
   2 \leq r<\infty \quad \text{and}\quad 1\leq p<2r/(r-1).
\end{equation*}
Let us denote that
\begin{equation*}
    D_{1}:= \dot{\mathcal{B}}^{2/p+1/r-2}_{p,q;r}(\mathbb{R}^3) 
    \quad \text{and} \quad 
    S_1:=\dot{\mathcal{B}}^{2/p+1/r-1}_{p,q;r}(\mathbb{R}^3).
\end{equation*}
Then the system \eqref{eq:NS2} is well-posed from $D_1$ to $\mathbb{P} S_1.$
\end{thm}

To extend the result of Theorem \ref{thm:wp_1} to the case $1\leq r<2,$ we need to bound some additional critical regularity of the solution as in the following theorem.
\begin{thm}\label{thm:wp_2} 
Let $1 \leq p_1 \leq p_2 \leq \infty,$ $1 \leq q_1\leq q_2 \leq \infty$ and $1\leq r_1<\infty$
satisfying one of the following conditions:
\begin{enumerate}
    \item $1\leq r_1 <2$ and $2/p_1+2/p_2+1/r_1>2,$
    \item $2\leq r_1<\infty$ and $1\leq p_1<2r_1/(r_1-1).$
\end{enumerate}
For such indices, let us denote that 
\begin{equation*}
    D_{2}:= \dot{\mathcal{B}}^{2/p_1+1/r_1-2}_{p_1,q_1;r_1}(\mathbb{R}^3)
    \quad \text{and} \quad 
    S_2:=\dot{\mathcal{B}}^{2/p_2-1}_{p_2,q_2;\infty}(\mathbb{R}^3) 
    \cap \dot{\mathcal{B}}^{2/p_1+1/r_1-1}_{p_1,q_1;r_1}(\mathbb{R}^3).
\end{equation*}
Then the system \eqref{eq:NS2} is well-posed from  $D_2$ to $\mathbb{P} S_2$.
\end{thm}

Let us give some comments on the well-posedness results obtained in Theorems \ref{thm:wp_1} and \ref{thm:wp_2}.
\begin{itemize}
    \item In Theorem \ref{thm:wp_1} or \ref{thm:wp_2}, we mainly prove the existence of the solution when the given force $\div f$ is small enough as the result \cite{Kan-Koz-Shi-19} in the standard homogenous Besov spaces setting. In particular, the solution $u$ of \eqref{eq:NS2} also satisfies
    \begin{equation}\label{es:wp}
        \|u\|_{S_k} \leq C \|f\|_{D_k}  
        \quad (k=1 \,\,\,\text{or}\,\,\, 2) 
    \end{equation}
    whenever $\|f\|_{D_k} $ is sufficiently small. 
    On the other hand, our results in Theorems \ref{thm:wp_1} and \ref{thm:wp_2} can be extended the general external force which is \emph{not} necessary in the divergence form.
\item In Theorems \ref{thm:wp_1} and \ref{thm:wp_2}, the solution of \eqref{eq:NS2} is constructed in the hybrid Besov space $\dot{\mathcal{B}}^{2/p+1/r-1}_{p,q;r}(\mathbb{R}^3)$ which is indeed critical for the problem \eqref{eq:NS2} with respect to the scaling structure \eqref{scaling} (see Remark \ref{rmk:space} below). Furthermore, our results in Theorems \ref{thm:wp_1} and \ref{thm:wp_2} is difficult to extend to the case of the two-dimensional flows as the critical regularity is too rough for the nonlinear terms in \eqref{eq:NS2}.
\end{itemize}

As we note that the space $\dot{\mathcal{B}}^{2/p+1/r-1}_{p,q;r}(\mathbb{R}^3)$ may contain some functions with singularities (see Remark \ref{rmk:space}), 
there is a natural question: whether the critical solution obtained in Theorem \ref{thm:wp_1} or \ref{thm:wp_2} becomes smooth if the smoothness of the given force $\div f$ is improved. Given better external forces as in the following theorem, we prove that the critical solution of \eqref{eq:NS2} not only becomes continuous in $\R^3$ but also decays at the infinity.
\begin{thm}\label{thm:decay}
Let $\delta\geq 0,$ $0<\sigma<1/2,$ $4/(3-2\sigma)<p\leq \infty$ and $p':=p/(p-1).$
Let us denote
\begin{equation*}
\begin{aligned}
 \|f\|_{D_{p,\sigma,\delta}}
 & :=  \n{f}_{\dot{\mathcal{B}}^{2/p+\delta}_{p,1;\infty}(\mathbb{R}^3)}^h
        +\n{\Jx^{\sigma} f}_{\dot{\mathcal{B}}^{2/p+1}_{p,1;\infty}(\mathbb{R}^3)}^h
        +\n{\Jx^{\sigma} f}_{\dot{\mathcal{B}}^{2/p-1-2\sigma}_{p,1;\infty}(\mathbb{R}^3)}^{\ell},\\
 \|u\|_{S_{p,\sigma,\delta}}
    &:= \n{u}_{\dot{\mathcal{B}}^{2/p+1+\delta}_{p,1;\infty}}^h
    +\n{\Jx^{\sigma} u }_{\dot{\mathcal{B}}^{2/p}_{p,1;\infty}}
\end{aligned}
\end{equation*}
with $\Jx:=(1+x_3^2)^{1/2}$ for any $x_3 \in \R.$ Suppose that 
\begin{equation*}
    \n{f}_{\dot{\mathcal{B}}^{2/p'-1}_{p',\infty;1}(\mathbb{R}^3)} 
   + \|f\|_{\dot{\mathcal{B}}^{2/p'-1-2\sigma}_{p',\infty;1}(\mathbb{R}^3)} 
   +\|f\|_{D_{p,\sigma,\delta}}
   \leq \eta
\end{equation*}
for some small positive constant $\eta.$ 
Then \eqref{eq:NS2} admits a unique solution $u$ satisfying 
\begin{equation}\label{des:0}
  \|u\|_{S_{p,\sigma,\delta}} \leq C  \|f\|_{D_{p,\sigma,\delta}}
\end{equation}
for some constant $C>0.$
\end{thm}

In particular, the estimate \eqref{des:0} and the embedding theory of the standard Besov spaces imply that
the solutions constructed by Theorem \ref{thm:decay} decay at the infinity with a speed of the polynomial with respect to the vertical variable. More precisely, we have 
\begin{equation*}
    |u(x_h,x_3)|
    \leq   
    \n{u(x_h,x_3)}_{L^{\infty}(\R^2_{x_h})}
    \lesssim  
    \n{u(x_h,x_3)}_{\dot{B}^{2/p}_{p,1}(\R^2_{x_h})} 
    \lesssim 
    \Jx^{-\sigma}
\end{equation*}
for $x_h=(x_1,x_2),$ $\Jx=(1+x_3^2)^{1/2}$ and $0<\sigma<1/2.$
Since the system \eqref{eq:NS} is isotropic, 
it is not hard to see the following result from Theorem \ref{thm:decay} by exchanging the order of $(x_1,x_2,x_3)$ and taking $\delta=0$. 
\begin{cor}
Let $0<\sigma<1/2,$ $4/(3-2\sigma)<p\leq \infty,$ $p':=p/(p-1)$ and set $\Jxx:= (1+|x|^2)^{1/2}$ for any $x\in \R^3.$
Suppose that 
\begin{equation*}
    \n{f}_{\dot{\mathbb{B}}^{2/p'-1}_{p',\infty;1}(\mathbb{R}^3)} 
   + \|f\|_{\dot{\mathbb{B}}^{2/p'-1-2\sigma}_{p',\infty;1}(\mathbb{R}^3)} 
   +\n{\Jxx^{\sigma} f}_{\dot{\mathbb{B}}^{2/p+1}_{p,1;\infty}(\mathbb{R}^3)}^h
        +\n{\Jxx^{\sigma} f}_{\dot{\mathbb{B}}^{2/p-1-2\sigma}_{p,1;\infty}(\mathbb{R}^3)}^{\ell}
   \leq \eta
\end{equation*}
for some small positive constant $\eta.$ 
Then \eqref{eq:NS2} admits a unique solution $u$ satisfying 
\begin{equation*}
|u(x)| \leq C \JX^{-\sigma}
\end{equation*}
for some constant $C>0.$
\end{cor}

Let us end up this section with outlining the structure of the rest part of this manuscript and  also the main idea in the proof. In Section \ref{sec:space}, we define the hybrid Besov space and prepare its important properties used in our proof. 
Next, in Section \ref{sec:linear}, we consider the linearized problem of \eqref{eq:NS2} where we transfer the corresponding linear operator to some ``evolution'' problem essentially characterized by some first-order semigroup.
Thus it is not hard to prove the boundedness of the linear operator in the hybrid Besov space. 
After that, we show the well-posedness in hybrid Besov spaces (Theorems 1.1 and 1.2) in Section \ref{sec:wp} by combining the results in Section \ref{sec:linear} and some product laws.
At last, we investigate the decay properties of the critical solution of \eqref{eq:NS2} in Section \ref{sec:decay}.
Inspired by the theory in \cite{Dan2000} for the non-stationary compressible Navier--Stokes flows, we decompose the solution into the low- and high-frequency parts. As in Lemma \ref{lemm:decay_Df_1}, the solution operator of the problem \eqref{eq:NS2} admits the different decay behaviors when we measure it in the low- and high-frequency norms.
In view of Lemma \ref{lemm:decay_Df_1} and some technical result (see Corollary \ref{cor:sc}), we can prove Theorem \ref{thm:decay} eventually.

\section{Functional spaces}
\label{sec:space}
In this section, we introduce the Besov type functional spaces which will be used in the whole paper. 
Let $\mathcal{S}(\R^N)$ and $\mathcal{S}'(\R^N)$ denote the space of all Schwartz functions and tempered distributions on $\R^N$ for $N=2,3.$
We set 
\begin{align*}
    \mathcal{S}_0(\R^N)
    :=
    \Mp{
    \phi \in \mathcal{S}(\R^N)
    \ :\ 
    \int_{\mathbb{R}^N}x^{\alpha}\phi(x)dx=0
    \ \text{for all }\alpha \in (\mathbb{N}\cup\{0\})^N
    }.
\end{align*}
It is well-known that the dual space $\mathcal{S}_0'(\R^N)$ of $\mathcal{S}_0(\R^N)$ is identified as 
\begin{align*}
    \mathcal{S}_0'(\R^N)  \sim
    \mathcal{S}'(\R^N)/\mathcal{P}(\mathbb{R}^N),
\end{align*}
where $\mathcal{P}(\mathbb{R}^N)$ denotes the set of all polynomials in $\R^N$.

Next, we recall the Littlewood--Paley theory. 
Let us consider some smooth radial functions $\varphi$ supported in the annulus $\{ \xi \in \R^2 : 3/4\leq |\xi| \leq {8}/{3} \}$ and satisfying the so-called \emph{Littlewood--Paley} decomposition
\begin{align*} 
 \sum_{j \in \mathbb{Z} } \varphi(2^{-j}\xi) =1,
 \quad \forall ~\xi \in \R^N \backslash \{0\}.
\end{align*}
Then we define  Fourier dyadic truncation operators as follows:
\begin{equation*}
\dot{\Delta}_j := \varphi(2^{-j}D), \quad  
\dot{S}_j := \sum_{j'\leq j-1} \dot{\Delta}_{j'},\,\,\forall j\in \Z.
\end{equation*}
Then the standard \emph{homogeneous} Besov norm of the distribution $f$ in the quotient space $\mathcal{S}'(\mathbb{R}^2)/ \mathcal{P}(\mathbb{R}^2)$ is given by 
\begin{equation*}
    \|f\|_{\dot{B}^{s}_{p,q}(\mathbb{R}^2)} 
    := \Big\|  \Big\{ 2^{js} \big\|\dot\Delta_j f(x) \big\|_{L^p_{x}} \Big \}_{j \in \mathbb{Z}} \Big\|_{\ell^q(\mathbb{Z})}
\end{equation*}
for any $s\in \R$ and $(p,q) \in [1,\infty]^2.$ 

The homogeneous Besov space is a very useful tool to analyze the Naiver-Stokes problem  (and other partial differential equations as well) in the whole space as it admits the following well-known scaling property:
\begin{equation}\label{scale:1}
    \|f(\lambda\cdot )\|_{\dot{B}^{s}_{p,q}(\mathbb{R}^2)} = \lambda^{s-2/p}  \|f(\cdot )\|_{\dot{B}^{s}_{p,q}(\mathbb{R}^2)} 
\end{equation}
for any constant $\lambda>0.$ Moreover, we shall use the following embedding property 
\begin{equation*}
    \dot{B}^{2/p}_{p,1}(\mathbb{R}^2) \hookrightarrow L^{\infty} (\R^2) \quad (1\leq p \leq \infty)
\end{equation*}
in the rest of this paper.
For more discussions on the theory of the standard Besov spaces, we refer to the monographs \cites{Bah-Che-Dan-11, Saw-18}.
\medskip

Based on the standard Besov spaces in $\R^2$, let us introduce the functional spaces with hybrid Besov regularities for the mappings defined in $\R^3:$ 
\begin{df}\label{def:Besov}
Let $s\in \R,$ and $(p,q,r) \in [1,\infty]^3.$ 
\begin{itemize}
    \item Let us set $x_h=(x_1,x_2)$ and define the Besov type space $\dot{\mathcal{B}}^{s}_{p,q;r}(\mathbb{R}^3)$ by the completion of $\mathcal{S}_0(\mathbb{R}^3)$ 
    with the norm 
    \begin{equation*}
    \|u\|_{\dot{\mathcal{B}}^{s}_{p,q;r}(\mathbb{R}^3)} 
    := \Big\|  \Big\{ 2^{js} \big\|\dot\Delta_j u(x_h,x_3) \big\|_{L^r_{x_3} L^p_{x_h} } \Big \}_{j \in \mathbb{Z}} \Big\|_{\ell^q(\mathbb{Z})}.
    \end{equation*}
    
    \item In a similar manner, we can also define the Besov type norm $\|\cdot\|_{ \dot{\mathcal{B}}^{s}_{p,q;r}(\mathbb{R}^3_{\sigma})} $ for the mappings 
    \begin{equation*}
        u : \mathbb{R}_{\sigma(1)} \rightarrow \mathcal{S}'\big( \mathbb{R}^2_{(x_{\sigma(2)},x_{\sigma(3)})} \big)/ 
        \mathcal{P}\big(\mathbb{R}^2_{(x_{\sigma(2)},x_{\sigma(3)})}\big)
    \end{equation*}
where $\sigma=\big(\sigma(1), \sigma(2),\sigma(3)\big)\in S_{\{1,2,3\}}$ 
stands for a permutation of $\{1,2,3\}.$
For simplicity, let us denote $\dot{\mathbb{B}}^{s}_{p,q;r}(\mathbb{R}^3)$ be the completion of $\mathcal{S}_0(\R^3)$ by the following norm
\begin{equation*}
    \|u\|_{ \dot{\mathbb{B}}^{s}_{p,q;r}(\mathbb{R}^3)} := \sum_{\sigma \in S_{\{1,2,3\}}} 
    \|u\|_{ \dot{\mathcal{B}}^{s}_{p,q;r}(\mathbb{R}^3_{\sigma})}<\infty.
\end{equation*}
\end{itemize}
\end{df}

\begin{rem}\label{rmk:space}
Let us give some comments on Definition \ref{def:Besov}.
\begin{itemize}
\item The space $\dot{\mathcal{B}}^{s}_{p,q;r}(\mathbb{R}^3)$ can be regarded as the extension of Chemin--Lerner spaces introduced in \cite{Che-Ler-95} to the steady problems. 
\item  The well-posedness results of \eqref{eq:NS2} in Theorems \ref{thm:wp_1} and \ref{thm:wp_2} are in the framework of the \emph{critical} functional spaces.
More precisely, in view of \eqref{scaling} and \eqref{scale:1},  the norm of $\dot{\mathcal{B}}^{2/p+1/r-1}_{p,q;r}(\mathbb{R}^3)$ is invariant under the transformation
    \begin{equation*}
        u(x)\leadsto  \lambda u( \lambda x), \quad \forall \,\,\,\lambda >0.
    \end{equation*}

\item For any fixed $r \in (1,\infty],$ there exists $\ep>0$ such that 
$$a_{\ep}:= 1-\ep-1/r>0.$$
Then the space $\dot{\mathcal{B}}^{2/p+1/r-1}_{p,1;r}(\mathbb{R}^3)$ may contain the functions with local singularity of the type $\chi (x_h,x_3)|x_h|^{-a_{\ep}},$
where $\chi$ is a suitable cut-off function near the origin.

\item As in \cite{Dan2000}, we decompose the function $u$ and its norm $\|u\|_{\dot{\mathcal{B}}^{s}_{p,q;r}(\mathbb{R}^3)}$ into the low- and high-frequency parts as follows:
\begin{equation*}
\begin{aligned}
    u^{\ell} (\cdot, x_3) &:= \sum_{j\leq N_0-1} \Dh_{j} u  (\cdot, x_3),\quad 
u^{h} (\cdot, x_3) := \sum_{j\geq N_0} \Dh_{j} u  (\cdot, x_3),\\
\|u\|_{\dot{\mathcal{B}}^{s}_{p,q;r}(\mathbb{R}^3)}^{\ell}&:=
 \Big\|  \Big\{ 2^{js} \big\|\dot\Delta_j u(x_h,x_3) \big\|_{L^r_{x_3} L^p_{x_h} } \Big \}_{j \leq N_0} \Big\|_{\ell^q(\mathbb{Z})},\\
 \|u\|_{\dot{\mathcal{B}}^{s}_{p,q;r}(\mathbb{R}^3)}^{h}&:=
 \Big\|  \Big\{ 2^{js} \big\|\dot\Delta_j u(x_h,x_3) \big\|_{L^r_{x_3} L^p_{x_h} } \Big \}_{j \geq N_0-1} \Big\|_{\ell^q(\mathbb{Z})}.
\end{aligned}
\end{equation*}
Above the integer $N_0$ is fixed. Analogously, we can also define the following notations:
\begin{equation*}
   \| \cdot \|_{{\dot{B}}^{s}_{p,q}(\mathbb{R}^2)}^h,\quad 
   \| \cdot \|_{{\dot{B}}^{s}_{p,q}(\mathbb{R}^2)}^{\ell},\quad 
     \| \cdot \|_{\dot{\mathbb{B}}^{s}_{p,q;r}(\mathbb{R}^3)}^h,\quad 
   \| \cdot \|_{\dot{\mathbb{B}}^{s}_{p,q;r}(\mathbb{R}^3)}^{\ell}.
\end{equation*}
For simplicity, we omit the detailed formulations.
\end{itemize}
\end{rem}

To handle the nonlinear terms in \eqref{eq:NS}, let us recall the definition of the paraproduct and remainder operators in horizontal variables as follows:
\begin{equation*}
\begin{aligned}
 \dot{T}_{f(\cdot, x_3)} g (\cdot, x_3) &:= \sum_{j \in \Z} 
    \dot{S}_{j-1} f (\cdot, x_3)\, \Dh_j g (\cdot, x_3), \\
    \dot{R} \big(f(\cdot, x_3),g(\cdot, x_3)\big) &:= \sum_{j\in \Z,\, |\nu|\leq 1} 
    \Dh_{j} f (\cdot, x_3) \,\Dh_{j+\nu} g (\cdot, x_3).
\end{aligned}
\end{equation*}
From Bony's calculus (see Chapter 2 of \cite{Bah-Che-Dan-11} for instance), we have the following technical result.
\begin{lemm}\label{lem:bony}
Let $k=1,2,$ $(s,s_1,s_2)\in \R^3$ and  $(p,q,r,p_k,q_k,r_k) \in [1,\infty]^6$  satisfying the following conditions:  
\begin{equation*}
    s=s_1+s_2, \quad
    1/p=1/p_1 +1/p_2, \quad 1/q=1/q_1 +1/q_2, \quad 1/r =1/r_1+1/r_2.
\end{equation*}
Then there exists a positive constant $C$ such that the following assertions hold. 
\begin{enumerate}
    \item For $s_1<0,$ we have 
    \begin{equation*}
    \begin{aligned}
       \big\| \dot{T}_{f} g  \big\|_{\dot{\mathcal{B}}^{s}_{p,q;r}(\mathbb{R}^3) } 
       &\leq C \|f\|_{\dot{\mathcal{B}}^{s_1}_{p_1,q_1;r_1}(\mathbb{R}^3) }
        \|g\|_{\dot{\mathcal{B}}^{s_2}_{p_2,q_2;r_2}(\mathbb{R}^3)},\\
       \big\| \dot{T}_{f} g  \big\|_{\dot{\mathcal{B}}^{s}_{p,q;r}(\mathbb{R}^3) } 
       &\leq C \|f\|_{L^{r_1}_{x_3}(L^{\infty}_{x_h}) }
        \|g\|_{\dot{\mathcal{B}}^{s}_{p,q;r_2}(\mathbb{R}^3)}.
    \end{aligned}
\end{equation*}
\item For $s_1+s_2>0,$ we have 
\begin{equation*}
    \big\| \dot{R}(f,g)  \big\|_{\dot{\mathcal{B}}^{s}_{p,q;r}(\mathbb{R}^3) }
    \leq C \|f\|_{\dot{\mathcal{B}}^{s_1}_{p_1,q_1;r_1}(\mathbb{R}^3) }
    \|g\|_{\dot{\mathcal{B}}^{s_2}_{p_2,q_2;r_2}(\mathbb{R}^3)}.
\end{equation*}
\end{enumerate}
\end{lemm}

\section{Solution of the linearized problem}
\label{sec:linear}
\subsection{Reformulation of the solution of \eqref{eq:NS2}}
According to the system \eqref{eq:NS2}, we begin with the following linearized problem 
\begin{align}\label{eq:Lin}
        -\Delta u = \mathbb{P}\div F
\end{align}
for some tensor $F=\{F_{jk}\}_{1\leq j,k\leq 3}.$
In particular, each component of $\eqref{eq:Lin}$ can be written as
\begin{equation}\label{eq:Lin-comp}
    -\Delta u_{\ell} = \sum_{j,k=1}^3 \left[ \delta_{\ell k} 
    +\partial_{x_\ell}\partial_{x_k}(-\Delta)^{-1} \right]
    \partial_{x_j}F_{kj} \quad (\ell=1,2,3).
\end{equation}
Using the equality 
\begin{equation}\label{eq:Laplace_operator}
    1 + \partial_{x_3}^2(-\Delta)^{-1} =-\Delta_h(-\Delta)^{-1}
\end{equation}
for $\Delta_h:=\partial_{x_1}^2+\partial_{x_2}^2,$
we obtain from \eqref{eq:Lin-comp} and \eqref{eq:Laplace_operator} that
\begin{align*}
    &\begin{aligned}
    -\Delta u_1
    ={}&
    \left[1+\partial_{x_1}^2(-\Delta)^{-1}\right]\sum_{j=1}^3 \partial_{x_j}F_{1j}
    +
    \partial_{x_1}\partial_{x_2}(-\Delta)^{-1}\sum_{j=1}^3\partial_{x_j}F_{2j}\\
    &
    +
    \partial_{x_1}\partial_{x_3}(-\Delta)^{-1}\sum_{j=1}^2\partial_{x_j}F_{3j}
    -
    \partial_{x_1}\left[1+\Delta_h(-\Delta)^{-1}\right]F_{33},
    \end{aligned}\\
    &\begin{aligned}
   -\Delta u_2
    ={}&
    \partial_{x_1}\partial_{x_2}(-\Delta)^{-1}\sum_{j=1}^3\partial_{x_j}F_{1j}
    +
    \left[1+\partial_{x_2}^2(-\Delta)^{-1}\right]\sum_{j=1}^3 \partial_{x_j}F_{2j}
    \\
    &
    +
    \partial_{x_2}\partial_{x_3}(-\Delta)^{-1}\sum_{j=1}^2\partial_{x_j}F_{3j}
    -
    \partial_{x_2}\left[1+\Delta_h(-\Delta)^{-1}\right]F_{33},
    \end{aligned}\\
    &\begin{aligned}
    -\Delta u_3
    ={}&
    \partial_{x_3}(-\Delta)^{-1}\sum_{j,k=1}^2\partial_{x_j}\partial_{x_k}F_{jk}
    -\sum_{j=1}^2\partial_{x_j}\left[1+\Delta_h(-\Delta)^{-1}\right]F_{j3}
    \\
    &
    -\Delta_h(-\Delta)^{-1}\sum_{j=1}^3\partial_{x_j}F_{3j}.
    \end{aligned}
\end{align*}
Then it is not hard to observe that all the terms in the formulation of $\mathbb{P}\div F$ are reduced to   
\begin{equation*}
\partial^{\alpha}g, \quad  (-\Delta)^{-1}\partial^{\beta}  g,
\end{equation*}
where the multi-indices  $\alpha=(\alpha_1,\alpha_2,\alpha_3)$
and $\beta=(\beta_h,\beta_3)=(\beta_1,\beta_2,\beta_3)
\in(\mathbb{N}\cup \{0\})^3$ satisfy 
$$|\alpha|=1, \quad |\beta_h|\geq 2 \quad \text{and}\quad |\beta|=3.$$ 
\medskip

Next, for convenience, let us introduce four model problems:
\begin{align}
    &
    -\Delta w^{(0)} 
    ={}
    g,\label{w0}\\
     &
    -\Delta w^{(1)} 
    ={}
    \partial_{x_3}g,\label{w1}\\
    &
    -\Delta \widetilde{w}^{(0)} 
    ={}
    (-\Delta)^{-1}g,\label{tw0}\\
    &
    -\Delta \widetilde{w}^{(1)} 
    ={}
    \partial_{x_3}(-\Delta)^{-1}g.\label{tw1}
\end{align}
Moreover, we denote the solution operators of  \eqref{w0}--\eqref{tw1} respectively by 
\begin{equation}\label{def:D_01}
   \mathcal{D}^{(0)}[g], \quad  \mathcal{D}^{(1)}[g], \quad 
   \widetilde{\mathcal{D}}^{(0)}[g], \quad \widetilde{\mathcal{D}}^{(1)}[g].
\end{equation}
\smallbreak 

In what follows, we derive the explicit formulas of the operators in \eqref{def:D_01}. 
For the operator $\mathcal{D}^{(0)},$  we apply the Fourier transform of \eqref{w0} with respect to $x_h$ to see that 
\begin{align}\label{Fw0}
    (|\xi_h|^2-\partial_{x_3}^2)\widehat{w^{(0)}}
    ={}
    \widehat{g}.
\end{align}
Since the corresponding eigen frequencies are $\pm|\xi_h|$, we can solve \eqref{Fw0} as
\begin{align*}
    \widehat{w^{(0)}}(\xi_h,x_3)
    ={}\widehat{\mathcal{D}^{(0)}}[g](\xi_h,x_3) 
    ={}\frac{1}{2|\xi_h|}
    \int_{\mathbb{R}}e^{-|x_3-y_3||\xi_h|}\widehat{g}(\xi_h,y_3)dy_3.
\end{align*}

For the solution operator $\mathcal{D}^{(1)}$ to \eqref{w1}, by the similar calculation for $\mathcal{D}^{(0)}$ and the integration by parts, we have
\begin{align*}
    \widehat{w^{(1)}}(\xi_h,x_3)
    ={} \widehat{\mathcal{D}^{(1)}}[g](\xi_h,x_3)    
    ={} -\frac{1}{2}
    \int_{\mathbb{R}}\operatorname{sgn}(x_3-y_3) e^{-|x_3-y_3||\xi_h|}\widehat{g}(\xi_h,y_3)dy_3.  
\end{align*}

For the solutions to \eqref{tw0} and \eqref{tw1}, notice that 
\begin{equation*}
    \widetilde{\mathcal{D}}^{(k)}= \mathcal{D}^{(0)} \circ \mathcal{D}^{(k)} \quad (k=0,1).
\end{equation*}
Then we combine the formulas of $\mathcal{D}^{(0)}$ and $\mathcal{D}^{(1)}$ above to obtain that 
\begin{align*}
    \widehat{\widetilde{w}^{(0)}}(\xi_h,x_3)
     ={}&
    \widehat{\widetilde{\mathcal{D}}^{(0)}}[g](\xi_h,x_3)\\
    ={}&
    \frac{1}{4|\xi_h|^2}
    \int_{\mathbb{R}}e^{-|x_3-y_3||\xi_h|}
    \int_{\mathbb{R}}e^{-|y_3-z_3||\xi_h|}
    \widehat{g}(\xi_h,z_3)dz_3dy_3,\\
    \widehat{\widetilde{w}^{(1)}}(\xi_h,x_3)
    ={}&
    \widehat{\widetilde{\mathcal{D}}^{(1)}}[g](\xi_h,x_3)\\
    ={}&-\frac{1}{4|\xi_h|}
    \int_{\mathbb{R}}e^{-|x_3-y_3||\xi_h|}
    \int_{\mathbb{R}}\operatorname{sgn}(y_3-z_3) e^{-|y_3-z_3||\xi_h|}\widehat{g}(\xi_h,z_3)dz_3dy_3.
\end{align*}

Now, taking advantage of the operators in \eqref{def:D_01}, the solution of the system \eqref{eq:Lin} can be formulated as follows:
\begin{align*}
    &\begin{aligned}
    \widehat{u_{\ell}}(\xi_h,x_3)
    ={}& \widehat{\mathcal{D}^{(0)}}
    \left[\sum_{j=1}^2 \partial_{x_j}F_{\ell j}-\partial_{x_\ell}F_{33}\right]
    +\widehat{\mathcal{D}^{(1)}}\left[F_{\ell 3}\right]\\
    &+\widehat{\widetilde{\mathcal{D}}^{(0)}}\left[
    \partial_{x_\ell}\left( \sum_{j,k=1}^2 \partial_{x_k}\partial_{x_j}F_{kj}
    -\Delta_h F_{33} \right)\right]
    +\widehat{\widetilde{\mathcal{D}}^{(1)}}\left[
    \partial_{x_\ell }\sum_{j=1}^2 \partial_{x_j}(F_{j3}+F_{3j}) \right]\\
    =:{}& \widehat{\mathcal{D}_\ell}[F](\xi_h,x_3), \quad \ell =1,2,
    \end{aligned}\\
    &\begin{aligned}
    \widehat{u_3}(\xi_h,x_3)
    ={}& 
    \widehat{\mathcal{D}^{(0)}}
    \left[
    -\sum_{j=1}^2\partial_{x_j}F_{j3}
    \right]
    +
    \widehat{\widetilde{\mathcal{D}}^{(0)}}
    \left[
    -\Delta_h\sum_{j=1}^2\partial_{x_j}(F_{j3}+F_{3j})
    \right]
    \\
    &
    +\widehat{\widetilde{\mathcal{D}}^{(1)}}\left[
    \sum_{j,k=1}^2\partial_{x_j}\partial_{x_k}F_{jk}
    -\Delta_h F_{33}
    \right]
    \\
    =:{}&
    \widehat{\mathcal{D}_3}[F](\xi_h,x_3).
    \end{aligned}
\end{align*}
For simplicity, we also write the solution $u=\mathcal{D}[F]$ of \eqref{eq:Lin} by setting 
\begin{equation}\label{def:D}
    \mathcal{D}[F] := \big( \mathcal{D}_1[F],\mathcal{D}_2[F],\mathcal{D}_3[F] \big)
\end{equation}
with $$\mathcal{D}_m[F]:=\mathscr{F}^{-1}_{x_h}\lp{\widehat{\mathcal{D}_m}[F]} 
\quad \text{for} \,\,\,m=1,2,3.$$

\begin{rem}\label{rem:def_D}
From the calculations above, the formula of the operator $\mathcal{D}$ can be regarded as the linear combination of the following operators
\begin{equation*}
     \mathcal{D}^{(0)} \nabla_h, \quad 
    \mathcal{D}^{(1)},\quad  
   \widetilde{\mathcal{D}}^{(0)}\nabla_h^{3}, \quad 
     \widetilde{\mathcal{D}}^{(1)}\nabla_h^{2}
\end{equation*}
with $\nabla_h=(\partial_1,\partial_2)^{\top}.$  
\end{rem}

Using notations above, we define the notion of the mild solution to \eqref{eq:NS2}.
\begin{df}
For a given external force $f$, we say that a vector field $u=(u_1(x),u_2(x),u_3(x))$ is a mild solution to \eqref{eq:NS} if $u$ satisfy 
\begin{align}\label{eq:mildNS}
    u = \mathcal{D}[f - u \otimes u].
\end{align}
\end{df}

\subsection{Boundedness of linear operators in \eqref{def:D_01}}
In this subsection, we derive the estimates of the linear operators introduced in \eqref{def:D_01}.
To this end, let us begin with the following technical lemma.
\begin{lemm}\label{lem:tech}
Let $1\leq p\leq \infty,$ $j\in \Z$ and $T>0.$
    There exist positive constants $c$ and $C$ independent of the choices of $T$ and $j$  such that the following estimates hold true.
    \begin{equation*}
        \begin{aligned}
            \n{|D_h|^{-1}\Dh_j f}_{L^p(\R^2_{x_h})} 
            & \leq  C 2^{-j}\n{\Dh_j f}_{L^p(\R^2_{x_h})},\\
      \n{e^{-|D_h|T}\Dh_j f}_{L^p(\R^2_{x_h})} 
            & \leq  Ce^{-cT2^{j}}\n{\Dh_j f}_{L^p(\R^2_{x_h})}.
        \end{aligned}
    \end{equation*}
Here, $x_h=(x_1,x_2)$  and
$g(|D_h|)$ denotes the Fourier multiplier in $\R^2_{x_h}$ with the symbol $g(|\xi_h|)$ for $g$ in $C^{\infty} (\R\backslash \{0\})$
and $\xi_h=(\xi_1,\xi_2).$
\end{lemm}
\begin{proof}
The first assertion is a direct consequence from Bernstein lemma (see Lemma 2.1 in \cite{Bah-Che-Dan-11} for example) while the second one was proved by \cite{Iwa-15}.
\end{proof}

Taking advantage of Lemma \ref{lem:tech}, 
we can establish the following results on the operators in \eqref{def:D_01}.
\begin{lemm} \label{lem:tech_2}
Let $j\in \Z$ and $(p,p_1,r,r_1,t_1) \in [1,\infty]^5$ satisfy
\begin{equation*}
    p_1\leq p \quad \text{and}\quad  1+1/r=1/r_1 +1/t_1.
\end{equation*}
There exists a positive constant $C$ such that the following assertions hold true.
\begin{enumerate}
    \item For $\mathcal{D}^{(k)}$ ($k=0,1$) defined by \eqref{def:D_01},  we have 
\begin{equation*}
        \n{ \Dh_j  \mathcal{D}^{(k)} [g] }_{L^r_{x_3}(L^p_{x_h})} 
        \leq C 2^{-(2-k+2/p-2/p_1+1/r-1/r_1)j} \n{ \Dh_j g }_{L^{r_1}_{x_3}(L^{p_1}_{x_h})}.
\end{equation*}
 \item In addition, let $(r_2,t_2) \in [1,\infty]^2$ satisfy
$$1+1/r_1=1/r_2 +1/t_2.$$
 For $\wt{\mathcal{D}}^{(k)}$ ($k=0,1$) defined by \eqref{def:D_01},  we have 
 \begin{equation*}
     \n{ \Dh_j  \widetilde{\mathcal{D}}^{(k)} [g] }_{L^r_{x_3}(L^p_{x_h})} 
     \leq C  2^{-(4-k+2/p-2/p_1+1/r-1/r_2)j}  \n{\Dh_j g }_{L^{r_2}_{x_3}(L^{p_1}_{x_h})}.
 \end{equation*}
\end{enumerate}
\end{lemm}

\begin{proof}
For convenience, we first derive the estimates of the operators $\mathcal{D}^{(0)}$ and $\widetilde{\mathcal{D}}^{(0)}$ by noticing the relationship 
\begin{equation}\label{eq:D_j-00}
    \Dh_j \widetilde{\mathcal{D}}^{(0)} [g]
=\Dh_j \mathcal{D}^{(0)}\mathcal{D}^{(0)} [g].
\end{equation}

Applying Lemma \ref{lem:tech}, Minkowski's inequality and Bernstein lemma to the following formulation
\begin{align*}
\Dh_j \mathcal{D}^{(0)} [g]
    =&\int_{x_3}^{\infty} \mathcal{F}_{\xi_h}^{-1} \left[   \frac{1}{2|\xi_h|}
    e^{-(y_3-x_3)|\xi_h|}\widehat{\Dh_j g}(\xi_h,y_3)  \right] (x_h) \,dy_3\\
    &+\int_{-\infty}^{x_3} \mathcal{F}_{\xi_h}^{-1} \left[   \frac{1}{2|\xi_h|}
    e^{-(x_3-y_3)|\xi_h|}\widehat{\Dh_j g}(\xi_h,y_3)  \right] (x_h) \,dy_3,
\end{align*}
we obtain that 
\begin{equation}\label{es:D_0}
\begin{aligned}
      \n{ \Dh_j  \mathcal{D}^{(0)} [g] (\cdot,x_3)}_{L^p_{x_h}} 
   &\lesssim  2^{-j} \int_{\R}   
    e^{-c|x_3-y_3|2^j} \|\Dh_j g (\cdot,y_3) \|_{L^p(\R^2)} \,dy_3\\
   & \lesssim  2^{-(1+2/p-2/p_1)j} \int_{\R}   
    e^{-c|x_3-y_3|2^j} \|\Dh_j g (\cdot,y_3) \|_{L^{p_1}(\R^2)} \,dy_3.
\end{aligned}
\end{equation}
Then  \eqref{es:D_0} and Young's inequality imply that 
\begin{equation}\label{es:D_0_1}
    \begin{aligned}
       \n{ \Dh_j  \mathcal{D}^{(0)} [g] }_{L^r_{x_3}(L^p_{x_h})} 
       & \lesssim 2^{-(1+2/p-2/p_1)j}   
       \n{e^{-c2^j y_3}}_{L^{t_1}_{y_3}(\R_+)}
                 \n{ \Dh_j g }_{L^{r_1}_{x_3}(L^{p_1}_{x_h})} \\
      &  \lesssim  2^{-(1+2/p-2/p_1+1/t_1)j} \n{ \Dh_j g }_{L^{r_1}_{x_3}(L^{p_1}_{x_h})} \\
      &  \lesssim  2^{-(2+2/p-2/p_1+1/r-1/r_1)j} \n{ \Dh_j g }_{L^{r_1}_{x_3}(L^{p_1}_{x_h})} 
    \end{aligned}
\end{equation}
for $1+1/r=1/r_1 +1/t_1$ and $\R_+:=(0,\infty).$
\smallbreak 

Next, from \eqref{eq:D_j-00} and \eqref{es:D_0_1}, we see that 
\begin{equation*}
    \begin{aligned}
        \n{ \Dh_j  \widetilde{\mathcal{D}}^{(0)} [g] }_{L^r_{x_3}(L^p_{x_h})} 
        &\lesssim 2^{-(2+2/p-2/p_1+1/r-1/r_1)j}
        \n{\Dh_j \mathcal{D}^{(0)} [g] }_{L^{r_1}_{x_3}(L^{p_1}_{x_h})} \\
        &\lesssim  2^{-(2+2/p-2/p_1+1/r-1/r_1)j} \cdot  2^{-(2+1/r_1-1/r_2)j}  
        \n{\Dh_j g }_{L^{r_2}_{x_3}(L^{p_1}_{x_h})} \\
        &\lesssim  2^{-(4+2/p-2/p_1+1/r-1/r_2)j}  \n{\Dh_j g }_{L^{r_2}_{x_3}(L^{p_1}_{x_h})} 
    \end{aligned}
\end{equation*}
with $1+1/r_1=1/r_2+1/t_2.$
\medskip 

At last, the estimates of $\mathcal{D}^{(1)}$ can be analogously established as the analysis of the operator $\mathcal{D}^{(0)}$ by using the formula
\begin{align*}
\Dh_j \mathcal{D}^{(1)} [g] 
    ={}& \frac{1}{2}
    \int_{x_3}^{\infty} 
    \mathcal{F}_{\xi_h}^{-1} \left[ e^{-(y_3-x_3)|\xi_h|}
    \widehat{\Dh_j  g}(\xi_h,y_3) \right](x_h)dy_3\\
     & - \frac{1}{2}
    \int_{-\infty}^{x_3} 
    \mathcal{F}_{\xi_h}^{-1} \left[ e^{-(x_3-y_3)|\xi_h|}
    \widehat{\Dh_j  g}(\xi_h,y_3) \right](x_h)dy_3.
\end{align*}
Indeed, it is not hard to see that 
\begin{equation}\label{es:D_1}
\begin{aligned}
      \n{ \Dh_j  \mathcal{D}^{(1)} [g] (\cdot,x_3)}_{L^p_{x_h}} 
   &\lesssim  \int_{\R}   
    e^{-c|x_3-y_3|2^j} \|\Dh_j g (\cdot,y_3) \|_{L^p(\R^2)} \,dy_3\\
   & \lesssim  2^{-j(2/p-2/p_1)} \int_{\R}   
    e^{-c|x_3-y_3|2^j} \|\Dh_j g (\cdot,y_3) \|_{L^{p_1}(\R^2)} \,dy_3,
\end{aligned}
\end{equation}
which implies that 
\begin{equation}\label{es:D_1_1}
    \begin{aligned}
       \n{ \Dh_j  \mathcal{D}^{(1)} [g] }_{L^r_{x_3}(L^p_{x_h})} 
    & \lesssim    2^{-(2/p-2/p_1)j} 
    \n{e^{-c2^j y_3}}_{L^{t_1}_{y_3}(\R_+)} 
      \n{ \Dh_j g }_{L^{r_1}_{x_3}(L^{p_1}_{x_h})} \\
      &  \lesssim  2^{-(2/p-2/p_1+1/t_1)j} \n{ \Dh_j g }_{L^{r_1}_{x_3}(L^{p_1}_{x_h})} \\
      &  \lesssim  2^{-(1+2/p-2/p_1+1/r-1/r_1)j} \n{ \Dh_j g }_{L^{r_1}_{x_3}(L^{p_1}_{x_h})} 
    \end{aligned}
\end{equation}
for $1+1/r=1/r_1 +1/t_1.$

On the other hand, we can obtain the bound of $\widetilde{\mathcal{D}}^{(1)}$ 
from \eqref{es:D_0_1}, \eqref{es:D_1_1} and the relationship 
$$\Dh_j \widetilde{\mathcal{D}}^{(1)} [g]
=\Dh_j \mathcal{D}^{(1)}\mathcal{D}^{(0)} [g].$$
For simplicity, we omit the details.
\end{proof}

Using Lemma \ref{lem:tech_2}, it is not hard to see the following result on the boundedness of the operator $\mathcal{D}$ defined by \eqref{def:D}.
\begin{lemm}\label{lemm:Df_1}
Let $s\in \R$ and $(p,p_1,q,q_1,r,r_1,t_1) \in [1,\infty]^7$ satisfy 
\begin{equation*}
    p_1\leq p,\quad q_1 \leq q \quad \text{and}\quad  1+1/r=1/r_1 +1/t_1.
\end{equation*}
There exists a positive constant $C$ such that  
\begin{equation}\label{es:Df_1}
    \|\mathcal{D}[F]\|_{\dot{\mathcal{B}}^{s}_{p,q;r}(\mathbb{R}^3)} 
    \leq C  \|F\|_{\dot{\mathcal{B}}^{s-1+\kappa}_{p_1,q_1 ;r_1}(\mathbb{R}^3)}
\end{equation}
with $\kappa :=2/p_1-2/p+1/r_1-1/r\geq 0.$
\end{lemm}

\begin{proof}
In order to establish \eqref{es:Df_1}, it suffices to prove that
\begin{equation}\label{es:lemm_Df_1}
     \|\Dh_j \mathcal{D} [g]\|_{L^r_{x_3}(L^p_{x_h})} 
    \lesssim 2^{-(1-\kappa)j}\|\Dh_j g\|_{L^{r_1}_{x_3}(L^{p_1}_{x_h})}
    \quad (j\in \Z),
\end{equation}
where $g$ stands for the general entry of the tensor $F.$ 
\smallbreak 

In fact, from Remark \ref{rem:def_D}, we only need to derive the $L^r_{x_3}(L^p_{x_h})$-bound of the following operators:
\begin{equation*}
   \Dh_j  \mathcal{D}^{(0)} \nabla_h, \quad 
   \Dh_j \mathcal{D}^{(1)},\quad  
    \Dh_j \widetilde{\mathcal{D}}^{(0)}\nabla_h^{3}, \quad 
    \Dh_j  \widetilde{\mathcal{D}}^{(1)}\nabla_h^{2}
\end{equation*}
with $\nabla_h=(\partial_1,\partial_2)^{\top}.$  
Then we apply Lemma \ref{lem:tech_2} and Bernstein lemma to obtain that 
\begin{equation*}
    \begin{aligned}
        \n{  \Dh_j  \mathcal{D}^{(0)} [\nabla_h g] }_{L^r_{x_3}(L^p_{x_h})}
     &\lesssim  2^{-(2-\kappa)j} \n{ \nabla_h \Dh_j g }_{L^{r_1}_{x_3}(L^{p_1}_{x_h})}
     \lesssim  2^{-(1-\kappa)j} \n{ \Dh_j g }_{L^{r_1}_{x_3}(L^{p_1}_{x_h})},\\
       \n{ \Dh_j  \mathcal{D}^{(1)} [g] }_{L^r_{x_3}(L^p_{x_h})} 
      &  \lesssim  2^{-(1-\kappa)j} \n{ \Dh_j g }_{L^{r_1}_{x_3}(L^{p_1}_{x_h})}, \\
            \n{ \Dh_j  \widetilde{\mathcal{D}}^{(0)}  
            [\nabla_h^{3} g] }_{L^r_{x_3}(L^p_{x_h})} 
   & \lesssim  2^{-(4-\kappa)j} \n{ \nabla_h^3 \Dh_j g }_{L^{r_1}_{x_3}(L^{p_1}_{x_h})}
   \lesssim  2^{-(1-\kappa)j}  \n{\Dh_j g }_{L^{r_1}_{x_3}(L^{p_1}_{x_h})},\\
     \n{ \Dh_j  \widetilde{\mathcal{D}}^{(1)}  
        [\nabla_h^{2} g] }_{L^r_{x_3}(L^p_{x_h})} 
   & \lesssim  2^{-(3-\kappa)j} \n{ \nabla_h^2 \Dh_j g }_{L^{r_1}_{x_3}(L^{p_1}_{x_h})}
   \lesssim  2^{-(1-\kappa)j}  \n{\Dh_j g }_{L^{r_1}_{x_3}(L^{p_1}_{x_h})}.
    \end{aligned}
\end{equation*}
This completes the proof of the estimate \eqref{es:lemm_Df_1}.
\end{proof}

\section{Well-posedness theory of \eqref{eq:NS2}}
\label{sec:wp}
In this section, we construct the mild solution of the problem \eqref{eq:NS2} in the sense of \eqref{eq:mildNS}.
In view of Lemma \ref{lemm:Df_1}, there hold that 
\begin{equation*}
 \|\mathcal{D}[f]\|_{S_k} \lesssim  \|f\|_{D_k} \quad  (k=1,2)
\end{equation*}
where $D_k$ and $S_k$ are introduced in Theorems \ref{thm:wp_1} and \ref{thm:wp_2} respectively.
Then to apply the standard well-posedness theory of the stationary problem as \cite{Kan-Koz-Shi-19},
the key point is to investigate the nonlinear part $\mathcal{D}[u\otimes u].$
In Subsections \ref{subsec:wp_1} and \ref{subsec:wp_2}, we derive the bound 
\begin{equation}\label{eq:non_r}
\n{\mathcal{D}[u\otimes u]}_{\dot{\mathcal{B}}^{2/p+1/r-1}_{p,q;r}(\mathbb{R}^3)}
\end{equation}
for $1\leq r<\infty,$ which will prove Theorems \ref{thm:wp_1} and \ref{thm:wp_2} respectively.
At last, in Subsection \ref{subsec:sc}, we end up with some theory on the propagation of some supercritical regularity which will be used later for the decay theory.

\subsection{Case $2 \leq r<\infty$}
\label{subsec:wp_1}
In this subsection, we consider \eqref{eq:non_r} for $2\leq r<\infty.$ 
\begin{lemm}\label{lem:nonlinear_1}
Let $(p,r) \in [1,\infty)^2$ fulfill the conditions 
\begin{equation}\label{cdt:pr_1}
   2 \leq r<\infty \quad \text{and}\quad 1\leq p<2r/(r-1). 
\end{equation}
 There exists a positive constant $C$ such that  
\begin{equation}\label{es:lem_uu_1}
     \|\mathcal{D}[u\otimes u]\|_{S_1} 
    \leq C  \n{u}_{S_1}^2
\end{equation}
for $S_1:=\dot{\mathcal{B}}^{2/p+1/r-1}_{p,q;r}(\mathbb{R}^3)$ as in Theorem \ref{thm:wp_1}.
\end{lemm}
\begin{proof}
According to Lemma \ref{lemm:Df_1}, we have 
\begin{equation}\label{es:lem_uu_2}
\|\mathcal{D}[u\otimes u]\|_{S_1}
    \lesssim   \n{u\otimes u}_{\dot{\mathcal{B}}^{2/p+2/r-2}_{p,q;r/2}(\mathbb{R}^3)}.
\end{equation}

Next, to bound the right hand side of \eqref{es:lem_uu_2}, we use the Bony's decomposition 
\begin{equation}\label{eq:bd}
    u_j u_k = \dot{T}_{u_j} u_k + \dot{T}_{u_k} u_j 
    +\dot{R}(u_j,u_k), \quad j,k=1,2,3.
\end{equation}
If there is no confusion, we denote $u_j$ and $u_k$ by $u$ for simplicity in what follows.
Then Lemma \ref{lem:bony} implies that 
\begin{equation*}
    \begin{aligned}
        \n{\dot{T}_u u}_{\dot{\mathcal{B}}^{2/p+2/r-2}_{p,q;r/2}(\mathbb{R}^3)}
        & \lesssim  \n{u}_{\dot{\mathcal{B}}^{1/r-1}_{\infty,\infty;r}(\mathbb{R}^3)}
        \n{u}_{\dot{\mathcal{B}}^{2/p+1/r-1}_{p,q;r}(\mathbb{R}^3)}
        \lesssim \n{u}_{S_1}^2.
    \end{aligned}
\end{equation*}

For the estimates of the remainder operator, we first consider the case $1\leq p \leq 2$ where $p'=p/(p-1)\geq p.$
According to Lemma \ref{lem:bony} and Bernstein lemma, we have 
\begin{equation*}
    \begin{aligned}
         \n{\dot{R}(u,u)}_{\dot{\mathcal{B}}^{2/p+2/r-2}_{p,q;r/2}(\mathbb{R}^3)}
         &\lesssim \n{\dot{R}(u,u)}_{\dot{\mathcal{B}}^{2/r}_{1,q;r/2}(\mathbb{R}^3)} \\
        & \lesssim  \n{u}_{\dot{\mathcal{B}}^{2/p'+1/r-1}_{p',\infty;r}(\mathbb{R}^3)}
        \n{u}_{\dot{\mathcal{B}}^{2/p+1/r-1}_{p,q;r}(\mathbb{R}^3)}
        \lesssim \n{u}_{S_1}^2. 
    \end{aligned}
\end{equation*}

On the other hand, let $2< p <2r/(r-1).$ Then Lemma \ref{lem:bony} implies that  
\begin{equation*}
    \begin{aligned}
         \n{\dot{R}(u,u)}_{\dot{\mathcal{B}}^{2/p+2/r-2}_{p,q;r/2}(\mathbb{R}^3)}
         &\lesssim \n{\dot{R}(u,u)}_{\dot{\mathcal{B}}^{4/p+2/r-2}_{p/2,q;r/2}(\mathbb{R}^3)} \\
        & \lesssim  \n{u}_{\dot{\mathcal{B}}^{2/p+1/r-1}_{p,\infty;r}(\mathbb{R}^3)}
        \n{u}_{\dot{\mathcal{B}}^{2/p+1/r-1}_{p,q;r}(\mathbb{R}^3)}
        \lesssim \n{u}_{S_1}^2.
    \end{aligned}
\end{equation*}
 Thus we conclude that 
\begin{equation}\label{es:lem_uu_3}
    \n{u\otimes u}_{\dot{\mathcal{B}}^{2/p+2/r-2}_{p,q;r/2}(\mathbb{R}^3)} 
    \lesssim  \n{u}_{S_1}^2
\end{equation}
for any indices $p$ and $r$ satisfying \eqref{cdt:pr_1}.
Therefore, we obtain \eqref{es:lem_uu_1} from \eqref{es:lem_uu_2} and \eqref{es:lem_uu_3}. This completes the proof.
\end{proof}

\subsection{Case $1\leq r<2$}
\label{subsec:wp_2}
For the propagation of the regularity in \eqref{eq:non_r} for $1\leq r<2,$  
it is difficult to use one single critical space as in Lemma \ref{lem:nonlinear_1}. 
Instead, we can prove the following result.
\begin{lemm}
Let $1 \leq p_1 \leq p_2 \leq \infty,$ $1 \leq q_1\leq q_2 \leq \infty$ and $1\leq r_1<\infty$
satisfy one of the following conditions:
\begin{enumerate}
    \item $1\leq r_1 <2$ and $2/p_1+2/p_2+1/r_1>2;$
    \item $2\leq r_1<\infty$ and $1\leq p_1<2r_1/(r_1-1).$
\end{enumerate}
Then there exists a constant $C$ such that
\begin{equation*}
\n{\mathcal{D}[u\otimes u]}_{S_2} \leq C \n{u}_{S_2}^2
\end{equation*}
for  $S_2:=\dot{\mathcal{B}}^{2/p_2-1}_{p_2,q_2;\infty}(\mathbb{R}^3) 
    \cap \dot{\mathcal{B}}^{2/p_1+1/r_1-1}_{p_1,q_1;r_1}(\mathbb{R}^3)$ as in Theorem \ref{thm:wp_2}.
\end{lemm}
\begin{proof}
According to Lemma \ref{lemm:Df_1}, we have 
\begin{equation*}
\begin{aligned} 
 \|\mathcal{D}[u\otimes u]\|_{S_2}   
 &\lesssim \begin{cases}
  \|u\otimes u\|_{\dot{\mathcal{B}}^{2/p_1+2/r_1-2}_{p_1,q_1;r_1/2}(\mathbb{R}^3)}
 & \text{if}\,\,\, 2 \leq r_1 <\infty;\\
 \|u\otimes u\|_{\dot{\mathcal{B}}^{2/p_1+1/r_1-2}_{p_1,q_1;r_1}(\mathbb{R}^3)}
 & \text{if}\,\,\, 1 \leq r_1 <2.
  \end{cases}
\end{aligned}
\end{equation*}
Moreover, if $r_1 \in [2,\infty)$ and $p_1 \in [1,2r_1/(r_1-1)),$ then \eqref{es:lem_uu_3} yields that 
\begin{equation*}
 \|\mathcal{D}[u\otimes u]\|_{S_2} 
 \lesssim \|u\otimes u\|_{\dot{\mathcal{B}}^{2/p_1+2/r_1-2}_{p_1,q_1;r_1/2}(\mathbb{R}^3)}
 \lesssim \|u\|_{\dot{\mathcal{B}}^{2/p_1+1/r_1-1}_{p_1,q_1;r_1}(\mathbb{R}^3)}^2
 \lesssim \n{u}_{S_2}^2.
\end{equation*}
\medskip

From now on, we assume that $1\leq r_1<2.$ 
Taking advantage of the convention below \eqref{eq:bd} and 
Lemma \ref{lem:bony}, we easily see that 
\begin{equation*}
    \begin{aligned}
        \n{\dot{T}_u u}_{\dot{\mathcal{B}}^{2/p_1+1/r_1-2}_{p_1,q_1;r_1}(\mathbb{R}^3)}
        & \lesssim  \n{u}_{\dot{\mathcal{B}}^{-1}_{\infty,\infty;\infty}(\mathbb{R}^3)}
        \n{u}_{\dot{\mathcal{B}}^{2/p_1+1/r_1-1}_{p_1,q_1;r_1}(\mathbb{R}^3)}\\
         & \lesssim  \n{u}_{\dot{\mathcal{B}}^{2/p_2-1}_{p_2,q_2;\infty}(\mathbb{R}^3)}
        \n{u}_{\dot{\mathcal{B}}^{2/p_1+1/r_1-1}_{p_1,q_1;r_1}(\mathbb{R}^3)}
        \lesssim \n{u}_{S_2}^2.
    \end{aligned}
\end{equation*}

If $1/p_1+ 1/p_2\geq 1,$ that is, $p_2 \leq p_1':=p_1/(p_1-1),$ then the estimate of the remainder term can be proved as follows 
\begin{equation*}
    \begin{aligned}
         \n{\dot{R}(u, u)}_{\dot{\mathcal{B}}^{2/p_1+1/r_1-2}_{p_1,q_1;r_1}(\mathbb{R}^3)}
        &\lesssim  \n{\dot{R}(u, u)}_{\dot{\mathcal{B}}^{1/r_1}_{1,q_1;r_1}(\mathbb{R}^3)}\\
           & \lesssim  \n{u}_{\dot{\mathcal{B}}^{2/p_1'-1}_{p_1',\infty;\infty}(\mathbb{R}^3)}
        \n{u}_{\dot{\mathcal{B}}^{2/p_1+1/r_1-1}_{p_1,q_1;r_1}(\mathbb{R}^3)}
        \lesssim \n{u}_{S_2}^2.
    \end{aligned}
\end{equation*}

On the other hand, suppose that $1/p_1+1/p_2=1/p_3<1$
for some $p_3 \in (1,\infty].$
Then Lemma \ref{lem:bony} implies that
\begin{equation*}
    \begin{aligned}
         \n{\dot{R}(u, u)}_{\dot{\mathcal{B}}^{2/p_1+1/r_1-2}_{p_1,q_1;r_1}(\mathbb{R}^3)}
        &\lesssim  \n{\dot{R}(u, u)}_{\dot{\mathcal{B}}^{2/p_3+1/r_1-2}_{p_3,q_1;r_1}(\mathbb{R}^3)}\\
           & \lesssim  \n{u}_{\dot{\mathcal{B}}^{2/p_2-1}_{p_2,\infty;\infty}(\mathbb{R}^3)}
        \n{u}_{\dot{\mathcal{B}}^{2/p_1+1/r_1-1}_{p_1,q_1;r_1}(\mathbb{R}^3)}
        \lesssim \n{u}_{S_2}^2
    \end{aligned}
\end{equation*}
so long as $2/p_3+1/r_1>2.$
\end{proof}

\subsection{Propagation of the supercritical regularity}
\label{subsec:sc}
In this subsection, we prove some technical result on the estimates of the \emph{supercritical} regularity compared the critical regularity in Theorems \ref{thm:wp_1} and \ref{thm:wp_2}.
\begin{prop} \label{prop:sc}
Let $0<\delta<1,$ $1 \leq p_1,p_2,q_1,q_2 \leq \infty$ and $1\leq r_1< \infty$
satisfy the following conditions
\begin{equation*}
p_1 \leq p_2,\quad q_1\leq q_2 \quad 
\text{and}\quad 2/p_2+2/p_1+(1-\delta)/r_1>2.
\end{equation*}
Suppose that $u \in \mathbb{P} \dot{\mathcal{B}}^{2/p_1+1/r_1-1}_{p_1,q_1;r_1}(\mathbb{R}^3)$ is a mild solution of \eqref{eq:mildNS} fulfilling the bound
\footnote{Note that it is possible to attain Condition \eqref{eq:c_f} by using Theorems \ref{thm:wp_1} and \ref{thm:wp_2}.}
\begin{equation}\label{eq:c_f}
    \n{u}_{\dot{\mathcal{B}}^{2/p_1+1/r_1-1}_{p_1,q_1;r_1}(\mathbb{R}^3)} \leq c_f.
\end{equation}
for some small constant $c_f$ depending on the choice of given force $f.$
Then there exists a constant $C$ depending on $c_f$ such that 
\begin{equation}\label{es:sc_1}
    \n{u}_{\dot{\mathcal{B}}^{2/p_2-1-\delta/r_1}_{p_2,q_2;\infty}(\mathbb{R}^3)}
    \leq C \|f\|_{\dot{\mathcal{B}}^{2/p_1-2+(1-\delta)/r_1}_{p_1,q_2;r_1}(\mathbb{R}^3)}
\end{equation}
provided, in addition, with $f\in \dot{\mathcal{B}}^{2/p_1-2+(1-\delta)/r_1}_{p_1,q_2;r_1}(\mathbb{R}^3).$
\end{prop}

\begin{proof}
Firstly, we observe from  Lemma \ref{lemm:Df_1} that 
    \begin{equation}\label{es:sc_2}
        \begin{aligned}
         \n{\mathcal{D}[f]}_{\dot{\mathcal{B}}^{2/p_2-1-\delta/r_1}_{p_2,q_2;\infty}(\mathbb{R}^3)}   
     & \lesssim \|f\|_{\dot{\mathcal{B}}^{2/p_1-2+(1-\delta)/r_1}_{p_1,q_2;r_1}(\mathbb{R}^3)},\\
     \n{\mathcal{D}[u\otimes u]}_{\dot{\mathcal{B}}^{2/p_2-1-\delta/r_1}_{p_2,q_2;\infty}(\mathbb{R}^3)} 
 & \lesssim \|u\otimes u\|_{\dot{\mathcal{B}}^{2/p_1-2+(1-\delta)/r_1}_{p_1,q_1;r_1}(\mathbb{R}^3)}.
        \end{aligned}
    \end{equation}

Next, to bound $u\otimes u$, we still apply Bony Calculus \eqref{eq:bd} as before.  
Then we easily see from Lemma \ref{lem:bony} and \eqref{eq:c_f} that 
\begin{equation*}
    \begin{aligned}
        \n{\dot{T}_u u}_{\dot{\mathcal{B}}^{2/p_1-2+(1-\delta)/r_1}_{p_1,q_1;r_1}(\mathbb{R}^3)}
        & \lesssim  \n{u}_{\dot{\mathcal{B}}^{-1-\delta/r_1}_{\infty,\infty;\infty}(\mathbb{R}^3)}
        \n{u}_{\dot{\mathcal{B}}^{2/p_1+1/r_1-1}_{p_1,q_1;r_1}(\mathbb{R}^3)}\\
         & \lesssim  c_f \n{u}_{\dot{\mathcal{B}}^{2/p_2-1-\delta/r_1}_{p_2,q_2;\infty}(\mathbb{R}^3)}.
    \end{aligned}
\end{equation*}

For the estimate of remainder operator, we consider the case $1/p_1+1/p_2 \geq 1$ where 
$p_2\leq p_1':=p_1/(p_1-1).$ Thus we have 
\begin{equation*}
    \begin{aligned}
         \n{\dot{R}(u, u)}_{\dot{\mathcal{B}}^{2/p_1-2+(1-\delta)/r_1}_{p_1,q_1;r_1}(\mathbb{R}^3)}
         & \lesssim \n{\dot{R}(u, u)}_{\dot{\mathcal{B}}^{(1-\delta)/r_1}_{1,q_1;r_1}(\mathbb{R}^3)}\\
           & \lesssim  \n{u}_{\dot{\mathcal{B}}^{2/p_1'-1-\delta/r_1}_{p_1',\infty;\infty}(\mathbb{R}^3)}
        \n{u}_{\dot{\mathcal{B}}^{2/p_1+1/r_1-1}_{p_1,q_1;r_1}(\mathbb{R}^3)}\\
         & \lesssim  c_f \n{u}_{\dot{\mathcal{B}}^{2/p_2-1-\delta/r_1}_{p_2,q_2;\infty}(\mathbb{R}^3)}.
    \end{aligned}
\end{equation*}

On the other hand, if $1/p_1+1/p_2 =1/p_3< 1$ for some $p_3 \in (1,\infty],$ then we have 
\begin{equation*}
    \begin{aligned}
         \n{\dot{R}(u, u)}_{\dot{\mathcal{B}}^{2/p_1-2+(1-\delta)/r_1}_{p_1,q_1;r_1}(\mathbb{R}^3)}
         & \lesssim \n{\dot{R}(u, u)}_{\dot{\mathcal{B}}^{2/p_3-2+(1-\delta)/r_1}_{p_3,q_1;r_1}(\mathbb{R}^3)}\\
           & \lesssim  \n{u}_{\dot{\mathcal{B}}^{2/p_2-1-\delta/r_1}_{p_2,\infty;\infty}(\mathbb{R}^3)}
        \n{u}_{\dot{\mathcal{B}}^{2/p_1+1/r_1-1}_{p_1,q_1;r_1}(\mathbb{R}^3)}\\
         & \lesssim  c_f \n{u}_{\dot{\mathcal{B}}^{2/p_2-1-\delta/r_1}_{p_2,q_2;\infty}(\mathbb{R}^3)}.
    \end{aligned}
\end{equation*}
so long as $2/p_3+(1-\delta)/r_1>2.$

Thus, combining all the bounds above, we obtain that 
\begin{equation}\label{es:sc_3}
    \|u\otimes u\|_{\dot{\mathcal{B}}^{2/p_1-2+(1-\delta)/r_1}_{p_1,q;r_1}(\mathbb{R}^3)} 
    \lesssim c_f  \n{u}_{\dot{\mathcal{B}}^{2/p_2-1-\delta/r_1}_{p_2,q;\infty}(\mathbb{R}^3)}.
\end{equation}
\medskip

At last, inserting \eqref{es:sc_3} to \eqref{es:sc_2}, we have
\begin{equation*}
\begin{aligned}
        \n{u}_{\dot{\mathcal{B}}^{2/p_2-1-\delta/r_1}_{p_2,q;\infty}(\mathbb{R}^3)}
    \leq C_1 \|f\|_{\dot{\mathcal{B}}^{2/p_1-2+(1-\delta)/r_1}_{p_1,q;r_1}(\mathbb{R}^3)}
    + C_2 c_f  \n{u}_{\dot{\mathcal{B}}^{2/p_2-1-\delta/r_1}_{p_2,q;\infty}(\mathbb{R}^3)}
\end{aligned}
\end{equation*}
for some constant $C_1$ and $C_2.$ 
Then taking $c_f$ small enough provides us with the bound \eqref{es:sc_1}. This completes our proof.
\end{proof}

In view of Theorem \ref{thm:wp_2} and Proposition \ref{prop:sc}, we have the following result.
\begin{cor} \label{cor:sc}
Let $0<\sigma<1,$ $1\leq p<4/(1+\sigma)$ and $1\leq q\leq \infty.$ Suppose that 
\begin{equation*}
   \n{f}_{\dot{\mathcal{B}}^{2/p-1}_{p,q;1}(\mathbb{R}^3)} 
   + \|f\|_{\dot{\mathcal{B}}^{2/p-1-\sigma}_{p,\infty;1}(\mathbb{R}^3)} 
   \leq \eta
\end{equation*}
for some small constant $\eta.$
Then there exists a unique solution $u$ of \eqref{eq:NS2} satisfying 
\begin{equation*}
    \n{u}_{\dot{\mathcal{B}}^{2/p-1}_{p,q;\infty}(\mathbb{R}^3)}
      + \n{u}_{\dot{\mathcal{B}}^{2/p}_{p,q;1}(\mathbb{R}^3)}+  \n{u}_{\dot{\mathcal{B}}^{2/p-1-\sigma}_{p,\infty;\infty}(\mathbb{R}^3)}
    \leq C \eta 
\end{equation*}
for some constant $C.$
\end{cor}
\begin{proof}
For any $f$ satisfying 
\begin{equation*}
    \n{f}_{\dot{\mathcal{B}}^{2/p-1}_{p,q;1}(\mathbb{R}^3)} \leq \eta \ll 1 \quad (1\leq p<4),
\end{equation*}
Theorem \ref{thm:wp_2} and \eqref{es:wp} imply that there exists a mild solution $u$ of \eqref{eq:NS2} fulfilling the estimate:
\begin{equation*}
    \n{u}_{\dot{\mathcal{B}}^{2/p-1}_{p,q;\infty}(\mathbb{R}^3)}
      + \n{u}_{\dot{\mathcal{B}}^{2/p}_{p,q;1}(\mathbb{R}^3)}  
    \lesssim  \eta.
\end{equation*}
Then we apply Proposition \ref{prop:sc} by taking 
$$\delta=\sigma,\quad p=p_1=p_2<4/(1+\sigma),\quad 
q_1=1,\quad q_2=\infty
\quad \text{and} \quad  r_1=1,$$ 
and we see that 
\begin{equation*}
    \n{u}_{\dot{\mathcal{B}}^{2/p-1-\sigma}_{p,\infty;\infty}(\mathbb{R}^3)}
    \lesssim  \|f\|_{\dot{\mathcal{B}}^{2/p-1-\sigma}_{p,\infty;1}(\mathbb{R}^3)} 
    \lesssim  \eta.
\end{equation*}

\end{proof}

\section{Pointwisely decay property of the solution of \eqref{eq:mildNS}}
\label{sec:decay}

This section is dedicated to the proof of Theorem \ref{thm:decay}.
To this end, we first derive some a priori bounds of $\mathcal{D}$ in the weighted norm in Subsection \ref{subsec:wn}, and then we study the solution of \eqref{eq:mildNS} by using the results obtained in Subsections \ref{subsec:sc}  and \ref{subsec:wn}.

\subsection{Some technical results on the weighted norms}
\label{subsec:wn}
Let us start with some bounds on the operators $\mathcal{D}^{(k)}$ and $\wt{\mathcal{D}}^{(k)}$ ($k=0,1$) given by \eqref{def:D_01}.
\begin{lemm}\label{lem:decay_tech_1}
Let $s\in \R,$ $\sigma>0$ and $(p,p_1,q,r,r_1,t_1) \in [1,\infty]^6$ satisfy 
\begin{equation*}
    p_1\leq p \quad \text{and}\quad  1+1/r=1/r_1 +1/t_1.
\end{equation*}
There exists a positive constant $C$ such that the following assertions hold true.
\begin{enumerate}
    \item Let us set $\kappa_1:=2/p_1-2/p+1/r_1-1/r\geq 0.$
    For $\mathcal{D}^{(k)}$ ($k=0,1$) defined by \eqref{def:D_01},  we have 
\begin{equation}\label{decay:D_k}
    \begin{aligned}
        \n{\Jx^{\sigma} \mathcal{D}^{(k)} [g] }_{\dot{\mathcal{B}}^{s}_{p,q;r}(\mathbb{R}^3)}^h 
        &\leq C \n{\Jx^{\sigma} g}_{\dot{\mathcal{B}}^{s+k-1+2/p_1-2/p}_{p_1,q;r_1}(\mathbb{R}^3)}^h,\\
   \n{\Jx^{\sigma} \mathcal{D}^{(k)} [g] }_{\dot{\mathcal{B}}^{s}_{p,q;r}(\mathbb{R}^3)}^{\ell} 
        &\leq C \n{\Jx^{\sigma} g}_{\dot{\mathcal{B}}^{s+k-2+\kappa_1-\sigma}_{p_1,q;r_1}(\mathbb{R}^3)}^{\ell}.
    \end{aligned}
\end{equation}

 \item In addition, let $(r_2,t_2) \in [1,\infty]^2$ satisfying 
$$1+1/r_1=1/r_2 +1/t_2.$$
Let us set $\kappa_2:=2/p_1-2/p+1/r_2-1/r\geq 0.$
 For $\wt{\mathcal{D}}^{(k)}$ ($k=0,1$) defined by \eqref{def:D_01},  we have 
\begin{equation*}
    \begin{aligned}
        \n{\Jx^{\sigma} \wt{\mathcal{D}}^{(k)} [g] }_{\dot{\mathcal{B}}^{s}_{p,q;r}(\mathbb{R}^3)}^h 
        &\leq C \n{\Jx^{\sigma} g}_{\dot{\mathcal{B}}^{s+k-2+2/p_1-2/p}_{p_1,q;r_2}(\mathbb{R}^3)}^h,\\
   \n{\Jx^{\sigma} \wt{\mathcal{D}}^{(k)} [g] }_{\dot{\mathcal{B}}^{s}_{p,q;r}(\mathbb{R}^3)}^{\ell} 
        &\leq C \n{\Jx^{\sigma} g}_{\dot{\mathcal{B}}^{s+k-4+\kappa_2-2\sigma}_{p_1,q;r_2}(\mathbb{R}^3)}^{\ell}.
    \end{aligned}
\end{equation*}
\end{enumerate}    
\end{lemm}

\begin{proof}
For any constant $\sigma>0$ and $x_3,y_3 \in \R,$  note that 
\begin{equation*}
    \Jx^\sigma \leq 2  \Jy^\sigma \Jxy^\sigma.
\end{equation*}
Then \eqref{es:D_0} and Young's inequality imply that 
\begin{multline}\label{des:D_0}
 \n{\Jx^{\sigma} \Dh_j  \mathcal{D}^{(0)} [g] (x_h,x_3)}_{L^r_{x_3}(L^p_{x_h})} 
    \lesssim  \, \n{2^{j\sigma} \Jx^{\sigma} e^{-c|x_3|2^j}}_{L^{t_1}_{x_3}(\R)} \\
    \times  2^{-(1+2/p-2/p_1+\sigma)j}  
    \n{\Jy^{\sigma}\Dh_j g (y_h,y_3) }_{L^{r_1}_{y_3}(L^{p_1}_{y_h})}
\end{multline}
for $1+1/r=1/r_1 +1/t_1.$

Suppose that $1\leq t_1<\infty.$ Then we notice that 
\begin{equation*}
    \begin{aligned}
    \int_{\R} 2^{j\sigma t_1} \Jx^{\sigma t_1} e^{-c|x_3|2^j t_1}  \,d x_3 
        = & \left(\int_{|x_3|\geq 1} +\int_{|x_3|<1} \right) 
        2^{j\sigma t_1} \Jx^{\sigma t_1} e^{-c|x_3|2^j t_1}  \,d x_3 \\
       \lesssim  &  2^{-j} t_1^{-1}e^{-c t_1 2^{j-1}} + 2^{j\sigma t_1},     
    \end{aligned}
\end{equation*}
which yields for any fixed $N_0\in \mathbb{N}$ that 
\begin{equation}\label{tech:hl}
 \n{2^{j\sigma} \Jx^{\sigma} e^{-c|x_3|2^j}}_{L^{t_1}_{x_3}(\R)}
       \lesssim 
      \begin{cases}
      2^{j\sigma} & \text{if}\,\,\,j\geq N_0-1,\\
    2^{-j(1+1/r-1/r_1)} & \text{if}\,\,\,j\leq  N_0.
      \end{cases}
\end{equation}
On the other hand, for the endpoint case where $t_1=r=\infty$ and $r_1=1,$ it is not hard to see that  
\begin{equation*}
    \begin{aligned}
    \n{2^{j\sigma} \Jx^{\sigma} e^{-c|x_3|2^j}}_{L^{\infty}_{x_3}(\R)}
    &\lesssim \sup_{|x_3|\geq 1}\big( (2^j |x_3|)^{\sigma} e^{-c|x_3|2^j} \big)
    +\sup_{|x_3|<1}\big( 2^{j\sigma} e^{-c|x_3|2^j} \big) \\
    & \lesssim  e^{-c 2^{j-1}} + 2^{j\sigma},
    \end{aligned}
\end{equation*}
which implies that \eqref{tech:hl} holds true for $t_1=\infty$ as well.

Thus \eqref{des:D_0} and \eqref{tech:hl} give us that 
\begin{equation}\label{decay:D_0_1}
\begin{aligned}
     &\n{\Jx^{\sigma} \Dh_j  \mathcal{D}^{(0)} [g] (x_h,x_3)}_{L^r_{x_3}(L^p_{x_h})} \\
    &\quad
    \lesssim        \begin{cases}
        2^{-(1+2/p-2/p_1)j}  \n{\Jy^{\sigma}\Dh_j g (y_h,y_3) }_{L^{r_1}_{y_3}(L^{p_1}_{y_h})}
        & \text{if}\,\,\,j\geq N_0-1,\\
    2^{-(2-\kappa_1+\sigma)j}   \n{\Jy^{\sigma}\Dh_j g (y_h,y_3) }_{L^{r_1}_{y_3}(L^{p_1}_{y_h})} & \text{if}\,\,\,j\leq  N_0,
      \end{cases}
\end{aligned}
\end{equation}
from which we prove \eqref{decay:D_k} for $k=0.$
\smallbreak 

For the operator $\mathcal{D}^{(1)},$ \eqref{es:D_1} and \eqref{tech:hl} yield that 
\begin{equation}\label{decay:D_1_1}
    \begin{aligned}
      & \n{\Jx^{\sigma} \Dh_j  \mathcal{D}^{(1)} [g] (x_h,x_3)}_{L^r_{x_3}(L^p_{x_h})} \\
    &\quad \lesssim  \n{2^{j\sigma} \Jx^{\sigma} e^{-c|x_3|2^j}}_{L^{t_1}_{x_3}(\R)} 
     2^{-(2/p-2/p_1+\sigma)j}  
    \n{\Jy^{\sigma}\Dh_j g (y_h,y_3) }_{L^{r_1}_{y_3}(L^{p_1}_{y_h})}\\
     & \quad   \lesssim          \begin{cases}
        2^{-(2/p-2/p_1)j}  \n{\Jy^{\sigma}\Dh_j g (y_h,y_3) }_{L^{r_1}_{y_3}(L^{p_1}_{y_h})}
        & \text{if}\,\,\,j\geq N_0-1,\\
    2^{-(1-\kappa_1+\sigma)j}   \n{\Jy^{\sigma}\Dh_j g (y_h,y_3) }_{L^{r_1}_{y_3}(L^{p_1}_{y_h})} & \text{if}\,\,\,j\leq  N_0.
      \end{cases}   
    \end{aligned}
\end{equation}
Then the bounds in \eqref{decay:D_k} also hold true for $k=1.$ 
\medskip

Next, to bound $ \widetilde{\mathcal{D}}^{(0)} [g]
=\mathcal{D}^{(0)}\mathcal{D}^{(0)} [g],$
it is not hard to see from \eqref{decay:D_0_1} that
\begin{equation*}
\begin{aligned}
     &\n{\Jx^{\sigma} \Dh_j  \wt{\mathcal{D}}^{(0)} [g] (x_h,x_3)}_{L^r_{x_3}(L^p_{x_h})} \\
    &\quad\lesssim          \begin{cases}
        2^{-(1+2/p-2/p_1)j}  \n{\Jy^{\sigma}\Dh_j \mathcal{D}^{(0)} g (y_h,y_3) }_{L^{r_1}_{y_3}(L^{p_1}_{y_h})}
        & \text{if}\,\,\,j\geq N_0-1,\\
    2^{-(2-\kappa_1+\sigma)j}   \n{\Jy^{\sigma}\Dh_j \mathcal{D}^{(0)}  g (y_h,y_3) }_{L^{r_1}_{y_3}(L^{p_1}_{y_h})} & \text{if}\,\,\,j\leq  N_0,
      \end{cases} \\
    &\quad \lesssim          \begin{cases}
        2^{-(2+2/p-2/p_1)j}  \n{\Jy^{\sigma}\Dh_j g (y_h,y_3) }_{L^{r_2}_{y_3}(L^{p_1}_{y_h})}
        & \text{if}\,\,\,j\geq N_0-1,\\
    2^{-(4-\kappa_2+2\sigma)j}   \n{\Jy^{\sigma}\Dh_j g (y_h,y_3) }_{L^{r_2}_{y_3}(L^{p_1}_{y_h})} & \text{if}\,\,\,j\leq  N_0.
      \end{cases}
\end{aligned}
\end{equation*}

In a similar manner, we can also obtain the bound of $ \widetilde{\mathcal{D}}^{(1)} [g]
=\mathcal{D}^{(1)}\mathcal{D}^{(0)} [g]$ from \eqref{decay:D_0_1} and \eqref{decay:D_1_1}. For simplicity, we omit the details.
This completes our proof.
\end{proof}

From Lemma \ref{lem:decay_tech_1}, it is not hard to see the following result.
\begin{lemm}\label{lemm:decay_Df_1}
Let $s\in \R,$ $\sigma >0$ and $(p,p_1,q,r,r_1,t_1) \in [1,\infty]^6$ satisfy
\begin{equation*}
    p_1\leq p \quad \text{and}\quad  1+1/r=1/r_1 +1/t_1.
\end{equation*}
Let us set $\kappa_1 :=2/p_1-2/p+1/r_1-1/r\geq 0.$
There exists a positive constant $C$ such that  
\begin{equation*}
\begin{aligned}
        \|\Jx^{\sigma}\mathcal{D}[F]\|_{\dot{\mathcal{B}}^{s}_{p,q;r}(\mathbb{R}^3)}^h 
    &\leq C  \|\Jx^{\sigma} F\|_{\dot{\mathcal{B}}^{s+1+2/p_1-2/p}_{p_1,q;r_1}(\mathbb{R}^3)}^h,\\
        \|\Jx^{\sigma}\mathcal{D}[F]\|_{\dot{\mathcal{B}}^{s}_{p,q;r}(\mathbb{R}^3)}^{\ell}
   & \leq C  \|\Jx^{\sigma} F\|_{\dot{\mathcal{B}}^{s-1+\kappa_1-2\sigma}_{p_1,q;r_1}(\mathbb{R}^3)}^{\ell}.     
\end{aligned}
\end{equation*}
\end{lemm}

\begin{proof}
Let $g$ denote the general entry of the tensor $F.$ 
Using Lemma \ref{lem:decay_tech_1} for $r_1=r_2$, we have 
\begin{equation*}
 \begin{aligned}
 \n{\Jx^{\sigma} \mathcal{D}^{(0)} [\nabla_h g] }_{\dot{\mathcal{B}}^{s}_{p,q;r}(\mathbb{R}^3)}^h 
 + \n{\Jx^{\sigma} \mathcal{D}^{(1)} [g] }_{\dot{\mathcal{B}}^{s}_{p,q;r}(\mathbb{R}^3)}^h 
  &\lesssim \n{\Jx^{\sigma} g}_{\dot{\mathcal{B}}^{s+2/p_1-2/p}_{p_1,q;r_1}(\mathbb{R}^3)}^h,\\
\n{\Jx^{\sigma} \wt{\mathcal{D}}^{(0)} [\nabla_h^{3}g] }_{\dot{\mathcal{B}}^{s}_{p,q;r}(\mathbb{R}^3)}^h 
 + \n{\Jx^{\sigma} \wt{\mathcal{D}}^{(1)} [\nabla_h^{2} g] }_{\dot{\mathcal{B}}^{s}_{p,q;r}(\mathbb{R}^3)}^h 
  &\lesssim \n{\Jx^{\sigma} g}_{\dot{\mathcal{B}}^{s+1+2/p_1-2/p}_{p_1,q;r_1}(\mathbb{R}^3)}^h,
    \end{aligned}
\end{equation*}
and 
\begin{equation*}
 \begin{aligned}
 \n{\Jx^{\sigma} \mathcal{D}^{(0)} [\nabla_h g] }_{\dot{\mathcal{B}}^{s}_{p,q;r}(\mathbb{R}^3)}^{\ell}
 + \n{\Jx^{\sigma} \mathcal{D}^{(1)} [g] }_{\dot{\mathcal{B}}^{s}_{p,q;r}(\mathbb{R}^3)}^{\ell} 
  &\lesssim \n{\Jx^{\sigma} g}_{\dot{\mathcal{B}}^{s-1+\kappa_1-\sigma}_{p_1,q;r_1}(\mathbb{R}^3)}^{\ell},\\
\n{\Jx^{\sigma} \wt{\mathcal{D}}^{(0)} [\nabla_h^{3}g] }_{\dot{\mathcal{B}}^{s}_{p,q;r}(\mathbb{R}^3)}^{\ell} 
 + \n{\Jx^{\sigma} \wt{\mathcal{D}}^{(1)} [\nabla_h^{2} g] }_{\dot{\mathcal{B}}^{s}_{p,q;r}(\mathbb{R}^3)}^h 
  &\lesssim \n{\Jx^{\sigma} g}_{\dot{\mathcal{B}}^{s-1+\kappa_1-2\sigma}_{p_1,q;r_1}(\mathbb{R}^3)}^{\ell}.
    \end{aligned}
\end{equation*}
Then it is not hard to obtain the desired bounds. 
\end{proof}

\subsection{Proof of Theorem \ref{thm:decay}}
Suppose that $0<\sigma <1/2$ and $4/(3-2\sigma)<p\leq \infty.$ 
Then, in view of Corollary \ref{cor:sc}, \eqref{eq:NS2} admits a solution $u$ satisfying 
\begin{equation}\label{Des:u_0}
    \n{u}_{\dot{\mathcal{B}}^{2/p'-1}_{p',\infty;\infty}(\mathbb{R}^3)}
      + \n{u}_{\dot{\mathcal{B}}^{2/p'}_{p',\infty; 1}(\mathbb{R}^3)}
      +\n{u}_{\dot{\mathcal{B}}^{2/p'-1-2\sigma}_{p,\infty;\infty}(\mathbb{R}^3)}
    \lesssim \eta
\end{equation}
so long as
\begin{equation*}
    \n{f}_{\dot{\mathcal{B}}^{2/p'-1}_{p',\infty;1}(\mathbb{R}^3)} 
   + \|f\|_{\dot{\mathcal{B}}^{2/p'-1-2\sigma}_{p',\infty;1}(\mathbb{R}^3)} 
   \leq \eta \ll 1.
\end{equation*}
Recall the notation that 
\begin{equation*}
\begin{aligned}
 \|u\|_{S_{p,\sigma,\delta}}
   &= \n{u}_{\dot{\mathcal{B}}^{2/p+1+\delta}_{p,1;\infty}}^h
    +\n{\Jx^{\sigma} u }_{\dot{\mathcal{B}}^{2/p}_{p,1;\infty}}\\ 
    &\sim  \n{u}_{\dot{\mathcal{B}}^{2/p+1+\delta}_{p,1;\infty}}^h
    +\n{\Jx^{\sigma} u }_{\dot{\mathcal{B}}^{2/p}_{p,1;\infty}}^h
    +\n{\Jx^{\sigma} u }_{\dot{\mathcal{B}}^{2/p}_{p,1;\infty}}^{\ell}
\end{aligned}
\end{equation*}
with $\delta \geq 0.$
In what follows, we shall derive the bound $\|u\|_{S_{p,\sigma,\delta}}$ for the solution $u$ fulfilling \eqref{Des:u_0}. For simplicity, let us denote $v:=\Jx^{\sigma} u$ in the rest of this subsection.
\medskip

Firstly, by \eqref{es:lemm_Df_1}, we easily see that 
\begin{equation*}
    \n{u}_{\dot{\mathcal{B}}^{2/p+1+\delta}_{p,1;\infty}(\mathbb{R}^3)}^h 
        \lesssim \n{f}_{\dot{\mathcal{B}}^{2/p+\delta}_{p,1;\infty}(\mathbb{R}^3)}^h
        + \n{u\otimes u }_{\dot{\mathcal{B}}^{2/p+\delta}_{p,1;\infty}(\mathbb{R}^3)}^h.
\end{equation*}
Then we use Lemma \ref{lem:bony} and \eqref{Des:u_0} to obtain that 
\begin{equation*}
\begin{aligned}
\n{u\otimes u }_{\dot{\mathcal{B}}^{2/p+\delta}_{p,1;\infty}}
&\lesssim \n{\dot T_u u}_{\dot{\mathcal{B}}^{2/p+\delta}_{p,1;\infty}} 
 + \n{\dot R(u,u)}_{\dot{\mathcal{B}}^{2/p+\delta}_{p,1;\infty}} \\
&\lesssim  \n{u}_{\dot{\mathcal{B}}^{-1}_{\infty,\infty;\infty}}
\n{u}_{\dot{\mathcal{B}}^{2/p+1+\delta}_{p,1;\infty}} \\
& \lesssim \n{u}_{\dot{\mathcal{B}}^{2/p'-1}_{p',\infty;\infty}(\mathbb{R}^3)} 
\big(\n{v}_{\dot{\mathcal{B}}^{2/p}_{p,1;\infty}}^\ell
+\n{u}_{\dot{\mathcal{B}}^{2/p+1+\delta}_{p,1;\infty}}^h\big)\\
&\lesssim \eta   \n{u}_{S_{p,\sigma,\delta}}.
\end{aligned}
\end{equation*}
Thus we have 
\begin{equation}\label{Des:u_1}
    \n{u}_{\dot{\mathcal{B}}^{2/p+1}_{p,1;\infty}(\mathbb{R}^3)}^h 
        \lesssim \n{f}_{\dot{\mathcal{B}}^{2/p+\delta}_{p,1;\infty}(\mathbb{R}^3)}^h
        + \eta   \n{u}_{S_{p,\sigma,\delta}}.
\end{equation}
In particular,  \eqref{Des:u_0} and \eqref{Des:u_1}  provide us with 
\begin{equation}\label{Des:u_2}
\begin{aligned}
   \n{u}_{\dot{\mathcal{B}}^{2/p+1}_{p,1;\infty}(\mathbb{R}^3)}
 &\lesssim \n{u}_{\dot{\mathcal{B}}^{2/p}_{p,1;\infty}(\mathbb{R}^3)}^{\ell} +  \n{u}_{\dot{\mathcal{B}}^{2/p+1+\delta}_{p,1;\infty}(\mathbb{R}^3)}^{h} \\
&\lesssim  \n{v}_{\dot{\mathcal{B}}^{2/p}_{p,1;\infty}(\mathbb{R}^3)}^{\ell}
+\n{f}_{\dot{\mathcal{B}}^{2/p+\delta}_{p,1;\infty}(\mathbb{R}^3)}^h
        + \eta   \n{u}_{S_{p,\sigma,\delta}}\\
&\lesssim  \n{f}_{\dot{\mathcal{B}}^{2/p+\delta}_{p,1;\infty}(\mathbb{R}^3)}^h
        + \n{u}_{S_{p,\sigma,\delta}}.
\end{aligned}
\end{equation}

Next, according to Lemma \ref{lemm:decay_Df_1}, $v$ admits the following a priori estimates :
\begin{equation}\label{des:1}
    \begin{aligned}
        \n{v}_{\dot{\mathcal{B}}^{2/p}_{p,1;\infty}(\mathbb{R}^3)}^h 
        & \lesssim \n{\Jx^{\sigma} f}_{\dot{\mathcal{B}}^{2/p+1}_{p,1;\infty}(\mathbb{R}^3)}^h
        + \n{v\otimes u }_{\dot{\mathcal{B}}^{2/p+1}_{p,1;\infty}(\mathbb{R}^3)}^h,\\
   \n{v}_{\dot{\mathcal{B}}^{2/p}_{p,1;\infty}(\mathbb{R}^3)}^{\ell} 
        & \lesssim \n{\Jx^{\sigma} f}_{\dot{\mathcal{B}}^{2/p-1-2\sigma}_{p,1;\infty}(\mathbb{R}^3)}^{\ell}
        +\n{v\otimes u}_{\dot{\mathcal{B}}^{2/p-1-2\sigma}_{p,1;\infty}(\mathbb{R}^3)}^{\ell}.
    \end{aligned}
\end{equation}

On one hand, by Minkowski's inequality, we observe that 
$S_{p,\sigma,\delta} \hookrightarrow L^{\infty}(\R^3)$
for any $p\in [1,\infty].$
Indeed, we have
\begin{equation}\label{eq:embed_1}
\n{(u,v)}_{L^{\infty}_{x_3}(L^{\infty}_{x_h})}
\lesssim \n{\n{v(\cdot,x_3)}_{\dot{\mathcal{B}}^{2/p}_{p,1}}}_{L^{\infty}_{x_3}}\lesssim \n{v}_{\dot{\mathcal{B}}^{2/p}_{p,1;\infty}}
\lesssim  \n{u}_{S_{p,\sigma,\delta}}.
\end{equation}
Then it is not hard to see from Lemma \ref{lem:bony}, \eqref{Des:u_2} and \eqref{eq:embed_1} that
\begin{equation}\label{nes:high}
\begin{aligned}
\n{v\otimes u }_{\dot{\mathcal{B}}^{2/p+1}_{p,1;\infty}}
&\lesssim \n{\dot T_v u}_{\dot{\mathcal{B}}^{2/p+1}_{p,1;\infty}} 
 + \n{\dot R(v,u)}_{\dot{\mathcal{B}}^{2/p+1}_{p,1;\infty}} \\
&\lesssim  \n{v}_{L^{\infty}_{x_3}(L^{\infty}_{x_h})}
\n{u}_{\dot{\mathcal{B}}^{2/p+1}_{p,1;\infty}}\\
&\lesssim  \n{f}_{\dot{\mathcal{B}}^{2/p+\delta}_{p,1;\infty}(\mathbb{R}^3)}^h
\n{u}_{S_{p,\sigma,\delta}} 
        +  \n{u}_{S_{p,\sigma,\delta}}^2 .
\end{aligned}
\end{equation}
Hence \eqref{des:1} and \eqref{nes:high} furnish that 
\begin{equation}\label{Des:v_1}
    \n{v}_{\dot{\mathcal{B}}^{2/p}_{p,1;\infty}(\mathbb{R}^3)}^h 
        \lesssim\n{\Jx^{\sigma} f}_{\dot{\mathcal{B}}^{2/p+1}_{p,1;\infty}(\mathbb{R}^3)}^h
 +\n{f}_{\dot{\mathcal{B}}^{2/p+\delta}_{p,1;\infty}(\mathbb{R}^3)}^h
\n{u}_{S_{p,\sigma,\delta}} 
        +  \n{u}_{S_{p,\sigma,\delta}}^2.
\end{equation}

On the other hand, note from Lemma \ref{lem:bony} that
\begin{equation*}
    \begin{aligned}
     \n{\dot T_u v}_{\dot{\mathcal{B}}^{2/p-1-2\sigma}_{p,1;\infty}} 
     &\lesssim  \n{u}_{\dot{\mathcal{B}}^{-1-2\sigma}_{\infty,\infty;\infty}}
    \n{v}_{\dot{\mathcal{B}}^{2/p}_{p,1;\infty}}
    \lesssim  \n{u}_{\dot{\mathcal{B}}^{2/p'-1-2\sigma}_{p',\infty;\infty}}
\n{v}_{\dot{\mathcal{B}}^{2/p}_{p,1;\infty}},\\
\n{\dot R(u,v)}_{\dot{\mathcal{B}}^{2/p-1-2\sigma}_{p,1;\infty}} 
&\lesssim \n{\dot R(u,v)}_{\dot{\mathcal{B}}^{1-2\sigma}_{1,1;1\infty}} 
\lesssim \n{u}_{\dot{\mathcal{B}}^{2/p'-1-2\sigma}_{p',\infty;\infty}}
    \n{v}_{\dot{\mathcal{B}}^{2/p}_{p,1;\infty}}
    \end{aligned}
\end{equation*}
since $0<\sigma<1/2.$
Then \eqref{Des:u_0} implies that 
\begin{equation}\label{nes:low}
\begin{aligned}
\n{v\otimes u}_{\dot{\mathcal{B}}^{2/{p}-1-2\sigma}_{p,1;\infty}}^{\ell}
&\lesssim 
\n{\dot T_u v}_{\dot{\mathcal{B}}^{2/p-1-2\sigma}_{p,1;\infty}} 
+\n{\dot R(u,v)}_{\dot{\mathcal{B}}^{2/p-1-2\sigma}_{p,1;1\infty}} \\
&\lesssim \n{u}_{\dot{\mathcal{B}}^{2/p'-1-2\sigma}_{p',\infty;\infty}}
\n{v}_{\dot{\mathcal{B}}^{2/p}_{p,1;\infty}} 
\lesssim \eta \n{u}_{S_{p,\sigma,\delta}}.
    \end{aligned}
\end{equation}
Thus we see from \eqref{des:1} and \eqref{nes:low} that 
\begin{equation}\label{Des:v_2}
    \begin{aligned}
    \n{v}_{\dot{\mathcal{B}}^{2/p}_{p,1;\infty}(\mathbb{R}^3)}^{\ell} 
        & \lesssim \n{\Jx^{\sigma} f}_{\dot{\mathcal{B}}^{2/p-1-2\sigma}_{p,1;\infty}(\mathbb{R}^3)}^{\ell}
        +\eta \n{u}_{S_{p,\sigma,\delta}}.       
    \end{aligned}
\end{equation}

At last, combining the bounds \eqref{Des:u_1}, \eqref{Des:v_1} and \eqref{Des:v_2} that 
\begin{equation*}
    \begin{aligned}
        \n{u}_{S_{p,\sigma,\delta}} 
        \lesssim & \n{f}_{\dot{\mathcal{B}}^{2/p+\delta}_{p,1;\infty}(\mathbb{R}^3)}^h
        +\n{\Jx^{\sigma} f}_{\dot{\mathcal{B}}^{2/p+1}_{p,1;\infty}(\mathbb{R}^3)}^h
        +\n{\Jx^{\sigma} f}_{\dot{\mathcal{B}}^{2/p-1-2\sigma}_{p,1;\infty}(\mathbb{R}^3)}^{\ell}\\
        &+ \big( \eta +  \n{f}_{\dot{\mathcal{B}}^{2/p+\delta}_{p,1;\infty}(\mathbb{R}^3)}^h\big) \n{u}_{S_{p,\sigma,\delta}}
        + \n{u}_{S_{p,\sigma,\delta}}^2.
    \end{aligned}
\end{equation*}
Then the bootstrap argument yields that 
\begin{equation*}
    \n{u}_{S_{p,\sigma,\delta}}
    \lesssim  \n{f}_{\dot{\mathcal{B}}^{2/p+\delta}_{p,1;\infty}(\mathbb{R}^3)}^h
        +\n{\Jx^{\sigma} f}_{\dot{\mathcal{B}}^{2/p+1}_{p,1;\infty}(\mathbb{R}^3)}^h
        +\n{\Jx^{\sigma} f}_{\dot{\mathcal{B}}^{2/p-1-2\sigma}_{p,1;\infty}(\mathbb{R}^3)}^{\ell}
\end{equation*}
by the smallness of the quantity $\|f\|_{D_{p,\sigma,\delta}}.$
This completes our proof of Theorem \ref{thm:decay}.

\section*{Acknowledgement}
MF is supported by Grant-in-Aid for Research Activity Start-up, Grant Number JP23K19011;
HT is supported by JSPS KAKENHI Grant Number JP24K16946;
XZ is partially supported by the National Natural Science Foundation of China (12101457) and the Fundamental Research Funds for the Central Universities.


\end{document}